\pdfoutput=1
\RequirePackage{ifpdf}
\ifpdf % We are running pdfTeX in pdf mode
\documentclass[pdftex]{sigma}
\else
\documentclass{sigma}
\fi

\numberwithin{equation}{section}
\newtheorem{Theorem}{Theorem}[section]
\newtheorem{Lemma}[Theorem]{Lemma}
\newtheorem{Proposition}[Theorem]{Proposition}
{\theoremstyle{definition}
\newtheorem{Definition}[Theorem]{Definition}
\newtheorem{Remark}[Theorem]{Remark}
}
\usepackage[all]{xy}

\begin{document}

\newcommand{\arXivNumber}{1302.0349}

\allowdisplaybreaks

\renewcommand{\thefootnote}{$\star$}

\renewcommand{\PaperNumber}{077}

\FirstPageHeading

\ShortArticleName{Quantitative~$K$-Theory Related to Spin Chern Numbers}

\ArticleName{Quantitative $\boldsymbol{K}$-Theory Related \\
to Spin Chern Numbers\footnote{This paper is a~contribution to the Special Issue on Noncommutative Geometry and Quantum
Groups in honor of Marc A.~Rief\/fel.
The full collection is available at
\href{http://www.emis.de/journals/SIGMA/Rieffel.html}{http://www.emis.de/journals/SIGMA/Rieffel.html}}}

\Author{Terry A.~LORING}

\AuthorNameForHeading{T.A.~Loring}

\Address{Department of Mathematics and Statistics, University of New Mexico,\\
Albuquerque, NM 87131, USA}
\Email{\href{mailto:loring@math.unm.edu}{loring@math.unm.edu}}
\URLaddress{\url{http://www.unm.edu/~loring/}}

\ArticleDates{Received January 15, 2014, in f\/inal form July 13, 2014; Published online July 19, 2014}

\Abstract{We examine the various indices def\/ined on pairs of almost commuting unitary matrices that can detect pairs
that are far from commuting pairs.
We do this in two symmetry classes, that of general unitary matrices and that of self-dual matrices, with an emphasis on
quantitative results.
We determine which values of the norm of the commutator guarantee that the indices are def\/ined, where they are equal,
and what quantitative results on the distance to a~pair with a~dif\/ferent index are possible.
We validate a~method of computing spin Chern numbers that was developed with Hastings and only conjectured to be
correct.
Specif\/ically, the Pfaf\/f\/ian--Bott index can be computed by the ``log method'' for commutator norms up to a~specif\/ic
constant.}


\Classification{19M05; 46L60; 46L80}

\begin{flushright}
\begin{minipage}{75mm}
\it Dedicated to Marc A.~Rieffel, whose lectures\\
 on Morita equivalence inspired all this
\end{minipage}
\end{flushright}

\renewcommand{\thefootnote}{\arabic{footnote}} \setcounter{footnote}{0}

\section{Introduction}

In the past decade, in condensed matter physics, certain systems with gapped Hamiltonians were found to fall into two
basic types.
Some were perturbations of completely trivial systems, and some were found to be far from all completely trivial
systems.
These are now called ``ordinary insulators'' and ``topological insulators'' respectively.
It is observed that a~path of systems perturbing a~topological insulator to an ordinary insulator must at some point
have closed the gap.
Physicists use~$K$-theory, both real and complex, to determine which insulators are which.

An older mathematical situation springs to mind here.
Given $C^{*}$-relations~\cite{LoringCstarRelations} in a~form where they can hold exactly and also hold approximately,
it was found that these approximate solutions fell into two basic types.
Some are close to exact solutions, and others are far away from all exact solutions.
The latter were called ``phantom approximate solutions''~\cite{EilersLoringContingenciesStableRelations}.
These approximate solutions were often found in matrix algebras $\mathbf{M}_{n}(\mathbb{C})$, but also in
$C^{*}$-algebras.
In either case, the main tools for distinguishing phantom from ordinary approximate solutions were constructions in
complex~$K$-theory.

The most basic set of $C^{*}$-relations are the relations, in the unital category,
\begin{gather*}
u^{*}u=1,\qquad uu^{*}=1,\qquad v^{*}v=1,\qquad vv^{*}=1,\qquad uv=vu.
\end{gather*}
There is generally no big distinction between~$u$ being ``almost unitary'' and being unitary, so we often study almost
commuting unitaries in $C^{*}$-algebras, meaning unitaries~$u$ and~$v$ with $  \Vert [u,v] \Vert
\leq\delta $ for some small~$\delta$ greater than zero.

There is a~direct connection between almost commuting unitary matrices and certain classes of f\/inite models of
topological insulators, explored in~\cite{HastLorTheoryPractice,LorHastHgTe,LorSorensenTorus}.
In that research, many interesting mathematical conjectures and questions were raised, some of which we address here.
Any serious numerical study of topological insulators must take into account the scattering
method~\cite{fulga2011scattering} of Fulga and his coauthors, which utilizes the sparseness of the matrices modeling
both position observables and the Hamiltonian.
That method still utilizes the Pfaf\/f\/ian--Bott index~\cite{LorHastHgTe}, discussed below, but as a~secondary calculation
after a~dimension reduction from 3D to 2D.

What we wish to emphasize here are aspects of phantom approximate solutions that are similar to the behavior in
topological insulators.
The connection between these two f\/ields of study is certainly greater than what has been explored to date.

Most of the theorems regarding approximate solutions to $C^*$-relations are completely non-quantitative.
There is often an constant $\delta_0$, unknown except for the fact that it is positive, so that nice things happen with
relations hold to within at most $\delta_0$.
One of the goals in this context is to develop ef\/f\/icient numerical algorithms.
When working numerically, a~constant like $10^{-10}$ can act ef\/fectively like zero.
Thus the desire for quantitative results.

A natural question regarding a~pair of unitary matrices that almost commute is: how close is this to a~pair that
actually commutes? The answer to this question necessarily involves the~$K$-theory of the two-torus.
Let us review how this connection arose.

There is a~particularly practical equation for the projection~$e$ in $\mathbf{M}_{2}(C(\mathbb{T}^{2}))$ that has rank
one and f\/irst Chern class one.
The formula is similar to that of the Rief\/fel projections~\cite{rieffel_A_theta} in the irrational rotation algebras,
specif\/ically
\begin{gather*}
e(z,w)=\left[
\begin{matrix}
 f(z) & g(z)+h(z)w
\\
g(z)+h(z)\overline{w} & 1-f(z)
\end{matrix}
\right],
\end{gather*}
where $f$, $g$ and~$h$ are certain real functions def\/ined on the unit circle.

The straight-forward plan in~\cite{loring1986torus} was to compute the~$K$-theory of a~$*$-homomorphism $\varphi:
C(\mathbb{T}^2) \rightarrow A$ by examining the associated commuting unitary elements $U=\varphi(u_0)$ and
$V=\varphi(v_0)$ of the AF algebra~$A$ and the projection
\begin{gather}
\label{eqn:basicBott}
e(U,V) = \left[
\begin{matrix} f(V) & g(z)+h(V)U
\\
g(V)+U^*h(V) & 1-f(V)
\end{matrix}
\right]
\end{gather}
where $u_0$ and $v_0$ are the canonical generating unitaries in $C(\mathbb{T}^2)$.
Unitaries in an AF algebra are well-known to be the limits of direct sums of unitary matrices, and so commuting unitaries are determined by sequences of almost commuting matrices.
In the specif\/ic situation of~\cite{loring1986torus}, $U=\lim (U_n \oplus A_n)$ and $V=\lim (V_n \oplus B_n)$ where
$A_n$, $B_n$ were commuting unitary matrices and $U_n$, $V_n$ were unitary matrices with $ \|[U_n, V_n
] \| \rightarrow 0$.

Equation~\eqref{eqn:basicBott} applies also to a~pair of almost commuting unitary matrices such as~$U_n$ and~$V_n$.
The result is not a~projection, but a~hermitian matrix with a~large gap at $\tfrac{1}{2}$ in its spectrum.
The~$K$-theory of~$\varphi$ was easily evaluated once the spectrum of $e(U_n,V_n)$ and $e(A_n,B_n)$ were understood.

The more interesting discovery in~\cite{loring1986torus} was that formula~\eqref{eqn:basicBott} can be used to def\/ine what
is now called the Bott index of a~pair of almost commuting unitary matrices.
This index can distinguish pairs of commuting matrices close to commuting pairs from those that are far from commuting
pairs.

There is ambiguity in the choice of~$f$,~$g$ and~$h$.
There are other ambiguities, discussed in~\cite{EilersLoringContingenciesStableRelations}, such that the fact that
$h(z)w$ could just as well been interpreted as
\begin{gather*}
\frac{1}{2} \{h(V),U \} =\frac{1}{2} (h(V)U+Uh(V) ).
\end{gather*}
To get good quantitative results about the distance to the closest commuting pair of unitary matrices, we will select
our functions and formulas very carefully.

In 1986 the only numerical computation of the Bott index that was practical involved relatively small matrices where~$V$
was diagonal.
Today we have from physics~\cite{fulga2011scattering,HastLorTheoryPractice,LorHastHgTe} large matrices where neither is
diagonal.
The cost of computing $f(V)$, $g(V)$ and $h(V)$ depends heavily on the choices in the scalar functions on the circle.

We end up with choices for~$f$,~$g$ and~$h$ that are very similar to the smooth functions illustrated
in~\cite{loring1986torus}, although we don't select them to have rapidly decreasing Fourier coef\/f\/icients.
We select functions that are well approximated by degree-f\/ive trigonometric polynomials and where the Fourier series are
relatively easy to calculate.

Soon after the Bott index was introduced, we found in joint work with Exel~\cite{ExelLoringInvariats} that a~simpler
formula based on winding numbers can be used.
We only proved that this formula worked for suf\/f\/iciently small commutator norms.
Here we will f\/ind a~concrete $\delta_{0}$ so that $\left\Vert [U,V]\right\Vert \leq\delta_{0}$ implies the
two invariants are equal.

We begin with a~survey, and some improved theorems, of the winding number index of~\cite{ExelLoringInvariats}.
Our expectation is that quantitative results regarding almost commuting matrices will be useful in applications,
especially in relation to topological insulators~\cite{fulga2011scattering,HastLorTheoryPractice,LorHastHgTe}.

We follow mathematical conventions, so $U^{*}$ refers to the conjugate-transpose.
To accommodate time reversal invariance in physics, we need to consider what in physics is called the dual operation,
\begin{gather}
\label{eq:dualDefn}
\left[
\begin{matrix} A & B
\\
C & D
\end{matrix}
\right]^{\sharp}=\left[
\begin{matrix} D^{\mathrm{T}} & -B^{\mathrm{T}}
\\
-C^{\mathrm{T}} & A^{\mathrm{T}}
\end{matrix}
\right].
\end{gather}

Specif\/ic unitary matrices that can be studied in the context of a~free particle system on a~f\/inite lattice on
a~two-torus are essentially complex-valued position operators that have been compressed to low energy space.
These are actually not quite unitary, but one can consider the unitary parts of their polar decomposition.
These then are almost commuting unitary matrices that carry a~lot of information about the original system~\cite[\S~1.1]{HastLorTheoryPractice}.

When the system has fermionic time reversal symmetry, the resulting unitary matrices will be self-dual.
The correct matrix problem to study is then almost commuting self-dual matrices.
The invariant~\cite{LorHastHgTe} that can be used to show that some pairs are bounded away from commuting self-dual
unitary pairs is the sign of
\begin{gather*}
\mathrm{Pf} (Q^{*}  (2e(U,V) - I ) Q  ),
\end{gather*}
where~$Q$ is a~specif\/ic matrix discussed below that creates anti-symmetry in the formula so that the Pfaf\/f\/ian makes
sense.
The resulting index we called the Pfaf\/f\/ian--Bott index~\cite{LorHastHgTe}.

The Pfaf\/f\/ian--Bott index of the unitary matrices associated to certain 2D systems has, for large system size, been
proven~\cite[Lemma~5.8.]{HastLorTheoryPractice} to equal the spin Chern number of that system.
It is perhaps more accurate to say that this Pfaf\/f\/ian--Bott index equals Kitaev's $\mathbb{Z}/2$ index for f\/inite 2D
systems in class~AII.

An alternate way to compute a~Bott index~\cite[Def\/inition~2.1]{ExelLoringInvariats}, equal to the Bott index for small
commutators, involves taking the logarithm of one of the unitaries.
It was surprising to f\/ind that this method, when adapted to the self-dual case, seemed to generate better data in
a~numerical study of disordered topological insulators~\cite{LorHastHgTe,LorSorensenTorus}.

We show in the f\/inal section that this method of computing the Pfaf\/f\/ian--Bott index gives the correct answer for commutator
norms up to a~specif\/ic constant.
There is a~separate issue of how to compute approximate logarithms of almost unitary matrices, and how to be sure to get
a~self-dual output given a~self-dual input.
That is discussed in a~separate paper~\cite{LoringLogOfUnitary}.

Our results are principally stated in terms of unitary matrices.
However, the study of almost commuting unitary elements of $C^{*}$-algebras is not that dif\/ferent.
We know this because we know that the soft-torus is RFD (residually f\/inite-dimensional)~\cite{EilersExelSoftTorusRFD}.

\section{The winding number invariant}

Given two unitary matrices~$U$ and~$V$ with $\delta= \Vert [U,V] \Vert $, we f\/ind
\begin{gather*}
 \Vert VUV^{*}U^{*}-I \Vert = \Vert  [U,V ] \Vert
\end{gather*}
and so by the spectral theorem
\begin{gather*}
\sigma (VUV^{*}U^{*} ) \subseteq  \{z\in\mathbb{T}\, |\,|z-1|\leq\delta  \}.
\end{gather*}
Thus when $\delta<2$ we can def\/ine $(VUV^{*}U^{*})^{t}$ for~$t$ between $0$ and~$1$, using a~branch of
$x^{t}$ with discontinuity on the negative~$x$-axis.
This is a~continuous path of unitary matrices from~$I$ to $VUV^{*}U^{*}$ and
\begin{gather*}
\det(VUV^{*}U^{*})=\det(I)=1,
\end{gather*}
so $ t\mapsto\det((VUV^{*}U^{*})^{t}) $ is a~loop on the unit circle.
We def\/ine $\omega(U,V)$ to be the winding number of this path.

This winding number invariant is very computable.
There are a~few alternate formulas, including
\begin{gather*}
\omega (U,V ) = \mathrm{Tr}\left(\frac{1}{2\pi i}\log (VUV^{*}U^{*} )\right)
%\label{eq:winding_via_log}
\end{gather*}
due to Exel~\cite[Lemma~3.1]{ExelSoftTorusI}.
This we easily prove: since in some basis $VUV^{*}U^{*}$ is diagonal and unitary,
\begin{gather*}
VUV^{*}U^{*}=\left(
\begin{matrix} e^{i\theta_{1}}
\\
& \ddots
\\
& & e^{i\theta_{n}}
\end{matrix}
\right)
\end{gather*}
for some $-\pi<\theta_{j}<\pi$ and
\begin{gather*}
\mathrm{Tr}\left(\frac{1}{2\pi i}\log\left(VUV^{*}U^{*}\right)\right) = \frac{1}{2\pi i}\mathrm{Tr} \left(
\begin{matrix} i\theta_{1}
\\
& \ddots
\\
& & i\theta_{n}
\end{matrix}
\right) = \frac{1}{2\pi}\sum\theta_{j}
\end{gather*}
and
\begin{gather*}
\det\left(\left(VUV^{*}U^{*}\right)^{t}\right) = \det\left(
\begin{matrix} e^{it\theta_{1}}
\\
& \ddots
\\
& & e^{it\theta_{n}}
\end{matrix}
\right) = \prod e^{it\theta_{j}}
\end{gather*}
and this path has winding number
\begin{gather*}
\omega\left(U,V\right)=\frac{1}{2\pi}\sum\theta_{j}.
\end{gather*}

\begin{Lemma}[\protect{\cite[p.~367]{ExelLoringInvariats}}] When~$U$ and~$V$ are commuting unitary matrices, $\omega(U,V)=0$.
\end{Lemma}

\begin{proof}
In this case the path of determinants is the constant path.
\end{proof}

\begin{Lemma}[\cite{ExelLoringAlmostCommutingUnitary}]
\label{lem:CyclicShift}
For
\begin{gather*}
U=\left(
\begin{matrix} 0 & & & & 1
\\
1 & 0
\\
& 1 & \ddots
\\
& & \ddots & 0
\\
& & & 1 & 0
\end{matrix}
\right),
\qquad
V=\left(
\begin{matrix} e^{i\pi/n}
\\
& e^{2i\pi/n}
\\
& & \ddots
\\
& & & e^{-i\pi/n}
\\
& & & & 1
\end{matrix}
\right)
\end{gather*}
we have $\omega(U,V)=-1$.
\end{Lemma}

\begin{proof}
We f\/ind $VUV^{*}U^{*}=e^{-\frac{2\pi i}{n}}I$ and so
\begin{gather*}
\frac{1}{2\pi i}\mathrm{Tr} (\log (VUV^{*}U^{*} ) ) = \frac{1}{2\pi i}\mathrm{Tr}\left(-\frac{2\pi
i}{n}I\right)=-1.\tag*{\qed}
\end{gather*}
\renewcommand{\qed}{}
\end{proof}

It is easy to modify the example in Lemma~\ref{lem:CyclicShift} to get a~pair of unitary matrices with $ \Vert
[U,V] \Vert =\delta$ and $\omega(U,V)=n$ for any $0<\delta<2$ and any~$n$.

\begin{Theorem}\label{thm:distance_to_commuting}
Consider a~pair of unitary matrices with $ \Vert [U,V] \Vert = \delta < 2$.
If $\omega(U,V)\neq 0$ then the distance to a~commuting pair of unitary matrices exceeds $\sqrt{2}$, meaning
\begin{gather*}
 \Vert U-U_{1} \Vert + \Vert V-V_{1} \Vert > \sqrt{2},
\end{gather*}
whenever $U_{1}$ and $V_{1}$ are unitary matrices with $U_{1}V_{1}=V_{1}U_{1}$.
Indeed,
\begin{gather*}
 \Vert U-U_{1} \Vert + \Vert V-V_{1} \Vert \geq \sqrt{2+\sqrt{4- \Vert
 [U,V ] \Vert^{2}}}.
\end{gather*}
\end{Theorem}

The proof of this will be broken into lemmas and propositions.
Theorem~\ref{thm:distance_to_commuting} is a~variation on the main result in~\cite{ExelLoringAlmostCommutingUnitary}.
That result had a~smaller lower bound, but the bound applied to the distance to any pair of commuting matrices, not just
commuting unitary matrices.

\begin{Lemma}
%\label{lem:gap_closing}
Suppose $
\Vert [U_{0},V_{0}]\Vert <2$.
If $\Vert [U_{1},V_{1}]\Vert <2$ and $ \omega(U_{0},V_{0})\neq\omega(U_{1},V_{1}) $ then for any
continuous path $U_s$ of unitary matrices from $U_0$ to $U_1$, and for any continuous path $V_s$ of unitary matrices
from $V_0$ to $V_1$, there must be at least one $s_{0}$ so that $\Vert [U_{s_{0}},V_{s_{0}}]\Vert
=2$.
\end{Lemma}

\begin{proof}
We will use a~homotopy argument.
If no such $s_{0}$ exists then
\begin{gather*}
(s,t)\mapsto\det \big( (V_{s}U_{s}^{*}V_{s}^{*}U_{s} )^{t} \big)
\end{gather*}
is a~homotopy between the path that determines $\omega(U_{0},V_{0})$ and the path that determines $\omega(U_{1},V_{1})$.
Therefore the winding numbers of these paths are equal.
\end{proof}

\begin{figure}[t]\centering
\includegraphics[clip,scale=0.5]{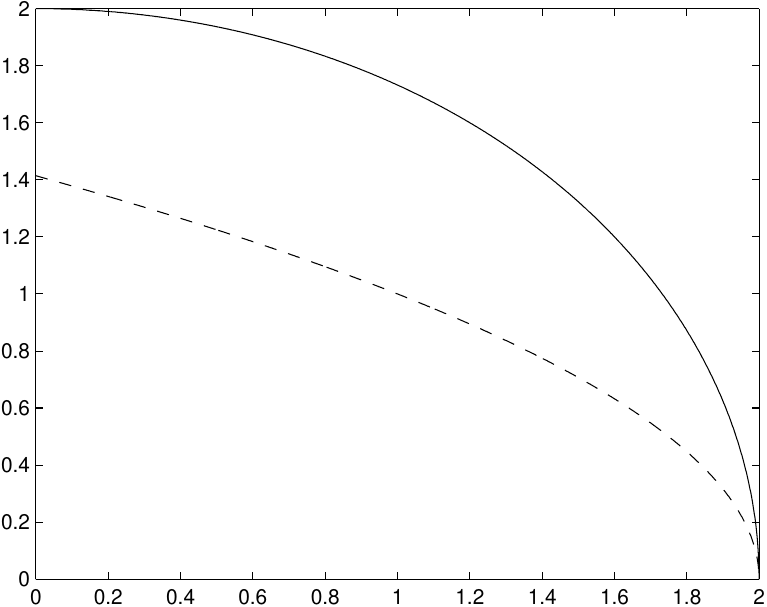}
\caption{The solid curve shows the minimum radius in the gap at $-1$
in with spectrum of $VUV^{*}U^{*}$ for unitaries~$U$ and~$V$, as a~function of $\delta= \Vert
[U,V] \Vert $. The lower curve show the minimum distance one must go to f\/ind a~pair with winding number index either undef\/ined or
dif\/ferent.}
\end{figure}

\begin{Proposition}
Suppose $ \Vert [U,V] \Vert <2$.
If $ \Vert [U_{1},V_{1}] \Vert =2$, then
\begin{gather*}
 \Vert U-U_{1} \Vert + \Vert V-V_{1} \Vert \geq \sqrt{2- \Vert  [U,V ] \Vert}.
\end{gather*}
If $ \Vert [U_{1},V_{1}] \Vert <2$ and $ \omega(U_{1},V_{1})\neq\omega(U,V), $ then
\begin{gather*}
 \Vert U-U_{1} \Vert + \Vert V-V_{1} \Vert \geq \sqrt{4-\big(\max( \Vert
 [U,V ] \Vert, \Vert  [U_{1},V_{1} ] \Vert)\big)^{2}}.
\end{gather*}
More generally, $ \Vert U-U_{1} \Vert + \Vert V-V_{1} \Vert $ is greater than or equal to
\begin{gather*}
\sqrt{2+2\sqrt{1-\frac{1}{4} \Vert  [U,V ] \Vert^{2}}\sqrt{1-\frac{1}{4} \Vert
 [U_{1},V_{1} ] \Vert^{2}}-\frac{1}{2} \Vert  [U,V ] \Vert  \Vert
 [U_{1},V_{1} ] \Vert}.
\end{gather*}
\end{Proposition}

\begin{proof}
We can connect $U=U_{0}$ to $U_{1}$ by an analytic path of unitary matrices $U_{t}$ of length
$2\arcsin(\frac{1}{2} \Vert U_{0}-U_{1} \Vert )$.
Similarly we have an analytic path $V_{t}$ from $V=V_{0}$ to $V_{1}$ of length $2\arcsin(\frac{1}{2} \Vert
V_{0}-V_{1} \Vert )$.
We now bound the length of the path $ W_{t}=V_{t}U_{t}V_{t}^{*}U_{t}^{*} $ in two ways.
We compute the derivative $W_{t}^{\prime}$ of $t \mapsto W_{t}$, and since $\Vert W_{t}^{\prime} \Vert \leq
2 \Vert U_{t}^{\prime} \Vert +2 \Vert V_{t}^{\prime} \Vert $ we see that
\begin{gather*}
\mathrm{Length}(W_{t}) \leq 4\arcsin\left(\frac{1}{2} \Vert U_{0}-U_{1} \Vert \right) +
4\arcsin\left(\frac{1}{2} \Vert V_{0}-V_{1} \Vert \right).
\end{gather*}
On the other hand, if $\mu_{k}(t)$ is an analytic choice of eigenvalues for $W_{t}$ then we know from~\cite[p.~374]{ExelLoringInvariats} that $ |\mu_{k}^{\prime}(t) |\leq \Vert W_{t}^{\prime} \Vert.
$ One of these paths of eigenvalues must hit~$-1$, and yet
\begin{gather*}
 |\mu_{k}(t)-1 |\leq \Vert  [U_{t},V_{t} ] \Vert
\end{gather*}
for $t=0$ and $t=1$, so
\begin{gather*}
\mathrm{Length}(W_{t})\geq2\pi-2\arcsin\left(\frac{1}{2} \Vert  [U_{0},V_{0} ] \Vert
\right)-2\arcsin\left(\frac{1}{2} \Vert [U_{1},V_{1} ] \Vert \right).
\end{gather*}
Therefore
\begin{gather*}
4\arcsin\frac{ \Vert U_{0}-U_{1} \Vert}{2} + 4\arcsin\frac{ \Vert V_{0}-V_{1} \Vert}{2} \geq 2\pi -
2\arcsin\frac{ \Vert  [U_{0},V_{0} ] \Vert}{2} - 2\arcsin\frac{ \Vert
 [U_{1},V_{1} ] \Vert}{2}.
\end{gather*}
We need to know the smallest value of $ \Vert U_{0}-U_{1} \Vert + \Vert V_{0}-V_{1} \Vert $ that can
be achieved, and for this it suf\/f\/ices to minimize these subject to the constraint
\begin{gather*}
4\arcsin\frac{ \Vert U_{0}-U_{1} \Vert}{2} + 4\arcsin\frac{ \Vert V_{0}-V_{1} \Vert}{2} = 2\pi -
2\arcsin\frac{ \Vert  [U_{0},V_{0} ] \Vert}{2} - 2\arcsin\frac{ \Vert
 [U_{1},V_{1} ] \Vert}{2}.
\end{gather*}
This is the problem of placing $6$ chords that are adjacent to each other that go around the unit circle, with two
chords f\/ixed of length $ \Vert [U_{0},V_{0}] \Vert $ and $ \Vert
[U_{1},V_{1}] \Vert $, while the other four come in pairs, two of length~$x$ and two of length~$y$.
The minimizing of $2x+2y$ occurs when we set one length, say~$y$, to zero, with the other the arc length corresponding
to arc length~$\pi$ minus half the arc length occupied by the two f\/ixed chords, so
\begin{gather*}
2\arcsin\left(\frac{1}{2}x\right) = \pi -\arcsin\left(\frac{1}{2} \Vert  [U_{0},V_{0} ] \Vert \right)
-\arcsin\left(\frac{1}{2} \Vert  [U_{1},V_{1} ] \Vert \right).
\end{gather*}
We conclude
\begin{gather*}
 \Vert U_{0}-U_{1} \Vert + \Vert V_{0}-V_{1} \Vert
\geq2\sin\left(\frac{1}{2}\left(\pi-\arcsin\left(\frac{1}{2} \Vert  [U_{0},V_{0} ] \Vert
\right)-\arcsin\left(\frac{1}{2} \Vert  [U_{1},V_{1} ] \Vert \right)\right)\right).
\end{gather*}
We f\/ind
\begin{gather*}
  2\sin\left(\frac{1}{2}\left(\pi-\arcsin\left(\frac{1}{2} \Vert  [U_{0},V_{0} ] \Vert
\right)-\arcsin\left(\frac{1}{2} \Vert  [U_{1},V_{1} ] \Vert \right)\right)\right)
\\
\qquad{}  =\sqrt{2}\sqrt{1-\cos\left(\pi-\arcsin\left(\frac{1}{2} \Vert  [U_{0},V_{0} ] \Vert
\right)-\arcsin\left(\frac{1}{2} \Vert [U_{1},V_{1} ] \Vert \right)\right)}
\\
 \qquad{}  =\sqrt{2}\sqrt{1+\cos\left(\arcsin\left(\frac{1}{2} \Vert  [U_{0},V_{0} ] \Vert
\right)+\arcsin\left(\frac{1}{2} \Vert  [U_{1},V_{1} ] \Vert \right)\right)}
\\
\qquad{}   =\sqrt{2+2\sqrt{1-\frac{1}{4} \Vert  [U_{0},V_{0} ] \Vert^{2}}\sqrt{1-\frac{1}{4} \Vert
 [U_{1},V_{1} ] \Vert^{2}}-\frac{1}{2} \Vert  [U_{0},V_{0} ] \Vert  \Vert
 [U_{1},V_{1} ] \Vert}
\end{gather*}
and so, setting $\Delta=\max( \Vert [U_{0},V_{0}] \Vert, \Vert
[U_{1},V_{1}] \Vert)$, we f\/ind
\begin{gather*}
 \Vert U_{0}-U_{1} \Vert + \Vert V_{0}-V_{1} \Vert
\geq\sqrt{2+2\sqrt{1-\frac{1}{4}\Delta^{2}}\sqrt{1-\frac{1}{4}\Delta^{2}}-\frac{1}{2}\Delta^{2}}
  =\sqrt{4-\Delta^{2}}.\tag*{\qed}
\end{gather*}
\renewcommand{\qed}{}
\end{proof}

We now get to a~dif\/f\/icult question.
Is the winding number invariant the only obstruction to closely approximating~$U$ and~$V$ by commuting unitary matrices?
It is important here that we stick with the operator norm in def\/ining ``close approximation'' as the answers to these
sort of questions can change dramatically if considering the Frobenius
norm~\cite{glebsky2010almost,gygi2003computation,marzariMLWF,RuheClosestNormal}.
(In particular, see the discussion in Section~III in~\cite{marzariMLWF}.) Results such as this also change dramatically
when the matrices come from dif\/ferent symmetry classes, as seen in~\cite{LorSorensenTorus}.

There is an answer, but it is only a~non-quantitative, nonconstructive result for small~$\delta$.
This we proven in joint work with Eilers and Pederson and we restate it here.
Also it matters that we are only interested in results for unitaries in $\mathbf{M}_{d}(\mathbb{C})$ that are
\emph{independent of~$d$}~\cite{HalmosUnknownDepthHilbert,LinAlmostCommutingMatrices}.

\begin{Theorem}[\protect{\cite[Theorem 6.15]{ELP-pushBusby}}]
%\label{thm:BottIsOnlyObstruction}
For any $\epsilon>0$, there is a~$\delta$ in $(0,2)$ so that, whenever~$U$ and~$V$ are unitary matrices in
$\mathbf{M}_{d}(\mathbb{C})$ with $ \Vert [U,V] \Vert \leq\delta$ and $\omega(U,V)=0$, there exist
unitary matrices $U_{1}$ and $V_{1}$ in $\mathbf{M}_{d}(\mathbb{C})$ so that
\begin{gather*}
 \Vert U-U_{1} \Vert + \Vert V-V_{1} \Vert \leq\epsilon
\end{gather*}
and $[U_{1},V_{1}]=0$.
\end{Theorem}

A serious limitation of the invariant $\omega(U,V)$ is that it does not generalize to unitaries in general
$C^*$-algebras, as it depends crucially on the determinant.
Another limitation is that we don't know how to modify it to work in other symmetry classes.
For example if we have self-dual unitary matrices, so $U^{\sharp}=U$ and $V^{\sharp}=V$, where $\sharp$ is a~specif\/ic
generalized involution detailed below, we f\/ind
\begin{gather*}
 (VUV^{*}U^{*} )^{\sharp}=U^{*}V^{*}UV
\end{gather*}
and so generally $VUV^{*}U^{*}$ is not self-dual.

\section[A direct~$K$-theory invariant~-- the Bott index]{A direct~$\boldsymbol{K}$-theory invariant~-- the Bott index}\label{sec:Bott-Index}

We need functional calculus of unitary matrices, also called matrix functions in applied mathematics.
An example is above where we applied the logarithm to a~unitary matrix.
Generally speaking, for the functional calculus $f(V)$ to be def\/ined for a~unitary matrix we need~$f$ def\/ined on the
circle.
One diagonalizes~$V$ via another unitary~$Q$ and applies~$f$ on the diagonal, so
\begin{gather*}
V=Q\left(
\begin{matrix} e^{i\theta_{1}}
\\
& \ddots
\\
& & e^{i\theta_{d}}
\end{matrix}
\right)Q^{*}\implies f(V)=Q\left(
\begin{matrix} f(e^{i\theta_{1}})
\\
& \ddots
\\
& & f(e^{i\theta_{2}})
\end{matrix}
\right)Q^{*}.
\end{gather*}
However, most of our calculations will involve Fourier series, and traditionally those are def\/ined in terms of scalar
functions that are periodic.

\begin{Definition}
Assume then that~$f$ is periodic of period $2\pi$ we def\/ine $f[V]$ as $\tilde{f}(V)$ where $\tilde{f}(z)=f(-i\log(z))$.
In other words,
\begin{gather*}
V=Q\left(
\begin{matrix} e^{i\theta_{1}}
\\
& \ddots
\\
& & e^{i\theta_{d}}
\end{matrix}
\right)Q^{*}\implies f[V]=Q\left(
\begin{matrix} f(\theta_{1})
\\
& \ddots
\\
& & f(\theta_{1})
\end{matrix}
\right)Q^{*}.
\end{gather*}
\end{Definition}

When~$f$ has uniformly convergent Fourier series, this is easier:
\begin{gather}
f(x)=\sum\limits_{n=-\infty}^{\infty}a_{n}e^{inx}\implies f[V]=\sum\limits_{n=-\infty}^{\infty}a_{n}V^{n}.
\label{eq:FourierSeriesFunctCalc}
\end{gather}

\begin{figure}[t]\centering
\includegraphics[clip,scale=0.45]{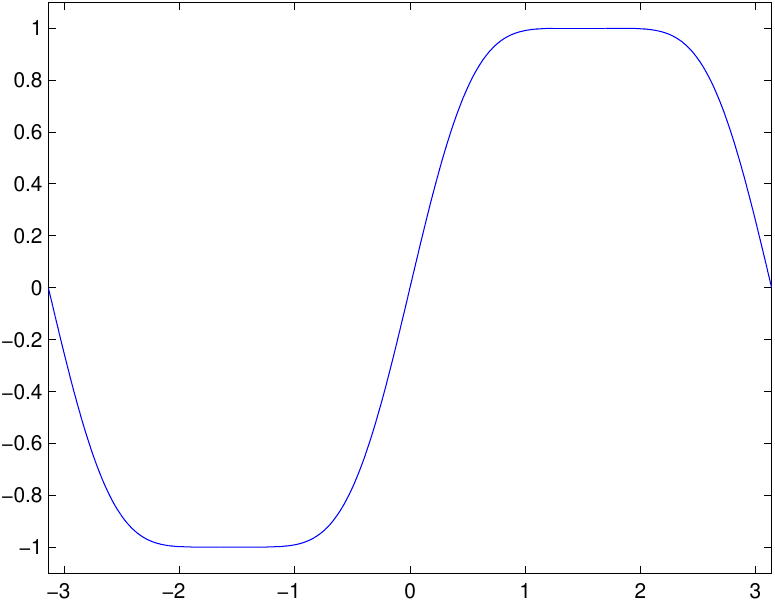}

\vspace{3mm}

\includegraphics[clip,scale=0.45]{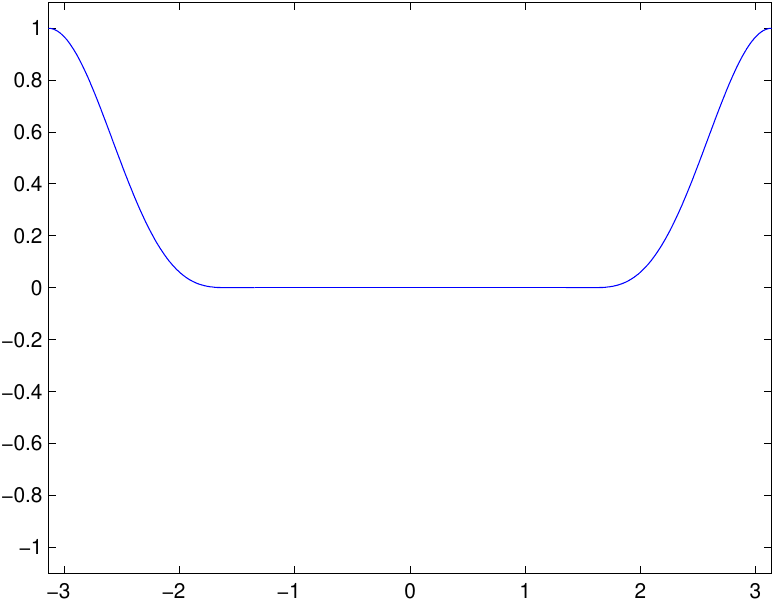}\qquad  \includegraphics[clip,scale=0.45]{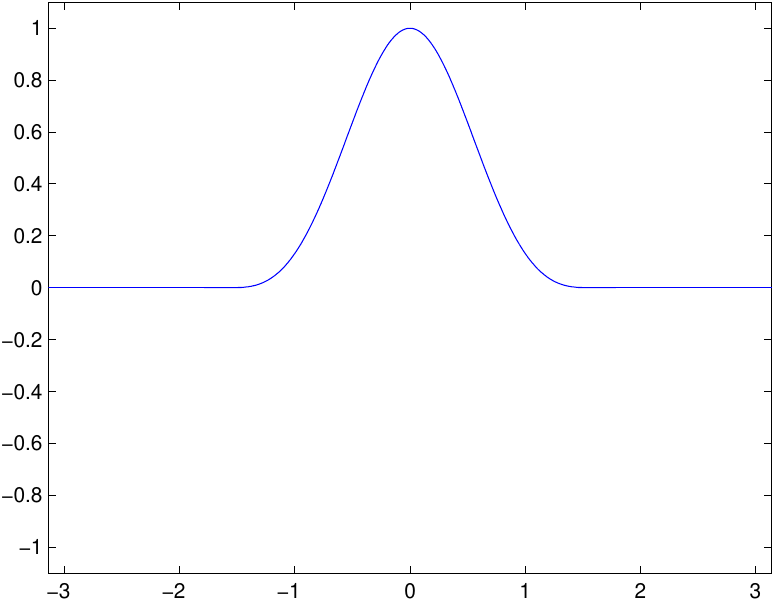}
\caption{Functions for the standard Bott index (trig method).}\label{fig:Functions-for-Bott}
\end{figure}

\begin{Definition}
\label{Def:BoffByf_g_h}
Def\/ine
\begin{gather*}
f(x)=\frac{1}{128}\left(150\sin(x)+25\sin(3x)+3\sin(5x)\right)
\end{gather*}
and
\begin{gather*}
g=
\begin{cases}
0, & x\in[-\frac{\pi}{2},\frac{\pi}{2}],
\\
\sqrt{1-f^{2}}, & x\notin[-\frac{\pi}{2},\frac{\pi}{2}],
\end{cases}
\qquad
h=
\begin{cases}
\sqrt{1-f^{2}}, & x\in[-\frac{\pi}{2},\frac{\pi}{2}],
\\
0, & x\notin[-\frac{\pi}{2},\frac{\pi}{2}],
\end{cases}
\end{gather*}
which are shown in Fig.~\ref{fig:Functions-for-Bott}.
For any unitaries set
\begin{gather*}
B(U,V)=\left(
\begin{matrix} f[V] & g[V]+\frac{1}{2}\left\{h[V],U\right\}
\vspace{1mm}\\
g[V]+\frac{1}{2}\left\{h[V],U^{*}\right\} & -f[V]
\end{matrix}
\right).
\end{gather*}
\end{Definition}

We have in mind unitary matrices, but let us look brief\/ly at the more abstract situation.
If we have commuting unitary matrices~$u$ and~$v$ in a~unital $C^*$-algebra~$A$ then we have, by the spectral theorem,
a~$*$-homomorphism $  \varphi: C(\mathbb T^2) \to B$.

Let us adopt the convention that $K_0$ will be def\/ined by hermitian elements with spectrum within $\{-1,1\}$ instead of
the usual description using projections, which are just hermitian elements with spectrum within $\{0,1\}$.
Then~$\varphi$ pushes forward the Bott element
\begin{gather*}
\beta = [B(z,w)] - \left[\left(
\begin{matrix} 1 & 0
\\
0 & -1
\end{matrix}
\right)\right]
\qquad \text{to}
\qquad
[B(u,v)] - \left[\left(
\begin{matrix} 1 & 0
\\
0 & -1
\end{matrix}
\right)\right].
\end{gather*}
If we have $\|[u,v]\|=\delta$ for small delta, then we can imagine something weaker than a~$*$-homomorphism, $
\psi: C(\mathbb T^2) \to B$, and attempt the push-forward.

Working heuristically, we simply create $[B(u,v)]$ and expect that this will be hermitian and with spectrum close to
being contained in $\{-1,1\}$.
Now we apply functional calculus $\chi(B(u,v)) $ (this is spectral f\/lattening in physics) and def\/ine the Bott index of
this pair as
\begin{gather*}
\chi(B(u,v)) - \left[\left(
\begin{matrix} 1 & 0
\\
0 & -1
\end{matrix}
\right)\right]
\end{gather*}
in $K_0(B)$.

This construction can be formalized in many ways.
Exel~\cite{ExelSoftTorusI} def\/ined the soft torus $A_\delta$ as the universal unital $C^*$-algebra generated by two
elements $u_\delta$ and $v_\delta$ subject to being unitary with $\|[u_\delta,v_\delta]\|\leq\delta$.
The only restriction on~$\delta$ is $\delta < 2$.
He calculated the~$K$-theory of $A_\delta$, showing that the natural map $\rho_\delta$ onto $C(\mathbb T^2) = A_0$ is an
isomorphism on~$K$-theory.
From~$u$ and~$v$ we get a~commuting diagram
\[
\xymatrix{
A_\delta \ar[rd]^(0.6){\gamma} \ar[d]^{\rho_\delta}\\
C(\mathbb T^2)  & B
}
\]
and can def\/ined a~very abstract index of $(u,v)$ as $\gamma_* \circ ((\rho_\delta)_*)^{-1} (\beta) $.

Computationally, spectrally f\/lattening an invertible matrix can be expensive.
Most importantly, doing so will destroy sparseness, should it initially exist.
Therefore, in the special case $B = \mathbf{M}_{n}(\mathbb{C})$ we use the signature.
Abstractly this is counting eigenvalues, but numerically there are many options for algorithms.

Now we resume discussion of the special case of unitary matrices.
For an invertible, hermitian matrix~$A$ we def\/ine its \emph{signature} $\mathrm{Sig}(A)$ as the number (with
multiplicity) of positive eigenvalues minus the number of negative eigenvalues.
We will prove that $\left\Vert [U,V]\right\Vert \leq0.206007$ forces $B(U,V)$ to be invertible.
Notice that if~$A$ is in $M_n(\mathbb C)$ and if~$n$ is even then $\mathrm{Sig}(A)$ must be even.

\begin{Remark}
Monte Carlo methods have generated numerical evidence that the gap in $B(U,V)$ actually closes at about $\delta = 0.85$.
\end{Remark}

\begin{Definition}
\label{defn:bott_index}
If $ \Vert [U,V] \Vert \leq0.206007$ def\/ine $\kappa(U,V)$, the {\em Bott index} of $(U,V)$, as the
integer
\begin{gather*}
\kappa(U,V)=\frac{1}{2}\mathrm{Sig} (B(U,V) ).
\end{gather*}
\end{Definition}

It should be noted that setting
\begin{gather*}
\gamma_{1}   =f(\theta_{2}),
\qquad
\gamma_{2}   =g(\theta_{2})+h(\theta_{2})\cos(\theta_{1}),
\qquad
\gamma_{3}   =h(\theta_{2})\sin(\theta_{1})
\end{gather*}
def\/ines the coordinates of a~map from $\mathbb{T}^{2}\rightarrow S^{2}\subseteq\mathbb{R}^{3}$ that has mapping degree
one.
Also notice $gh=0$ and $f^{2}+g^{2}+h^{2}=1$.

\begin{figure}[t]\centering
\includegraphics[clip,scale=0.55]{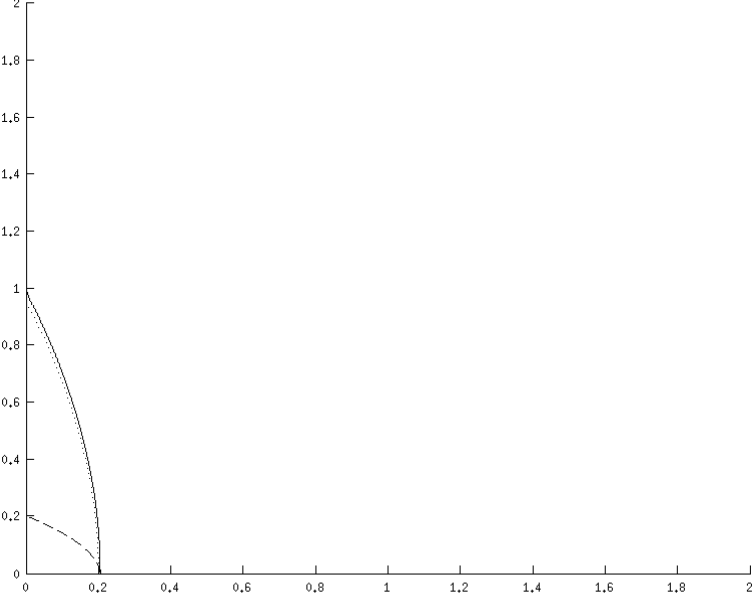}\qquad \includegraphics[clip,scale=0.55]{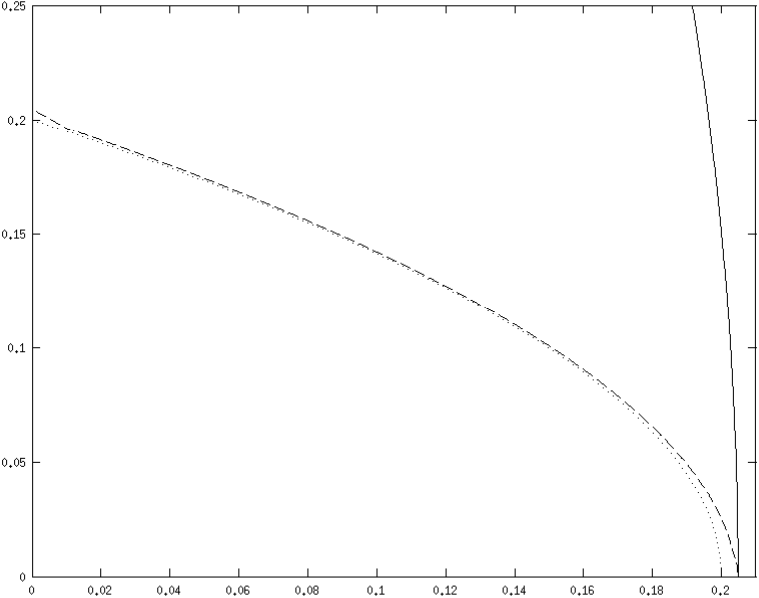}

\caption{In both: solid curve on top is the guaranteed spectral gap in $B(U,V)$ for a~given~$\delta$.
The dashed curve is the distance one can move where it is proven the gap will not close.
The dotted curve in the left plot is $\frac{19}{20}\sqrt{1-5\delta}$.
The dotted curve on the right is  $\frac{1}{5}\sqrt{1-5\delta}$.}\label{fig:smallGap}
\end{figure}

\begin{Theorem}
\label{thm:measure_the_gap}
Suppose~$U$ and~$V$ are unitaries and $ \delta= \Vert [U,V] \Vert \leq0.206007$.
\begin{enumerate}\itemsep=0pt
\item[$1.$] The hermitian matrix $B(U,V)$ has a~spectral gap at $0$ of radius at least $\frac{19}{20}\sqrt{1-5\delta}$.
Indeed, the gap is at least as large as the function of~$\delta$ plotted as a~solid curve in Fig.~{\rm \ref{fig:smallGap}}.
\item[$2.$] The distance $ \Vert U-U_{1} \Vert + \Vert V-V_{1} \Vert $ needed so that $B(U_{1},V_{1})$ has
$0$ in its spectrum is at least $\frac{1}{5}\sqrt{1-5\delta}$.
Indeed, this distance is at least as large as the function of~$\delta$ plotted as a~dashed curve in
Fig.~{\rm \ref{fig:smallGap}}.
\end{enumerate}
\end{Theorem}

We will prove Theorem~\ref{thm:measure_the_gap} in Section~\ref{sec:proof_of_gap}, and it will be a~lot of work.
Moreover, the gap here is much smaller than we saw for $VUV^{*}U^{*}$.
Why do we bother? The point is symmetry.

Suppose $U^{\sharp}=U$ and $V^{\sharp}=V$ for unitaries in $\mathbf{M}_{2d}(\mathbb{C})$ and with $\sharp$ the
generalized involution discussed in the next section, that physicists call the dual.
Then $B(U,V)$ is in $\mathbf{M}_{2d}(\mathbb{C})\otimes\mathbf{M}_{2}(\mathbb{C})$ which has on it the generalized
involution $\tau=\sharp\otimes\sharp$.
In terms of real $C^{*}$-algebras, this is a~copy of $\mathbf{M}_{2d+2}(\mathbb{C})$ with the transpose operation.
In physics language, we are tensoring two half-odd-integer spin systems to get a~system with integer spin, in
a~non-standard basis.
We f\/ind $ B(U,V)^{\tau}=-B(U,V) $ and so $B(U,V)$ def\/ines a~class in
\begin{gather*}
K_{2}(\mathbb{R})\cong K_{-2}(\mathbb{H})\cong\mathbb{Z}/2
\end{gather*}
that is computed directly in terms of the Pfaf\/f\/ian, hence the Pfaf\/f\/ian--Bott index studied in~\cite{LorHastHgTe}.
In return for a~small gap, indeed no guaranteed gap if $ \Vert [U,V] \Vert $ is too large, we get
a~construction that is amenable to symmetries.

\begin{figure}[t]\centering
\includegraphics[clip,scale=0.6]{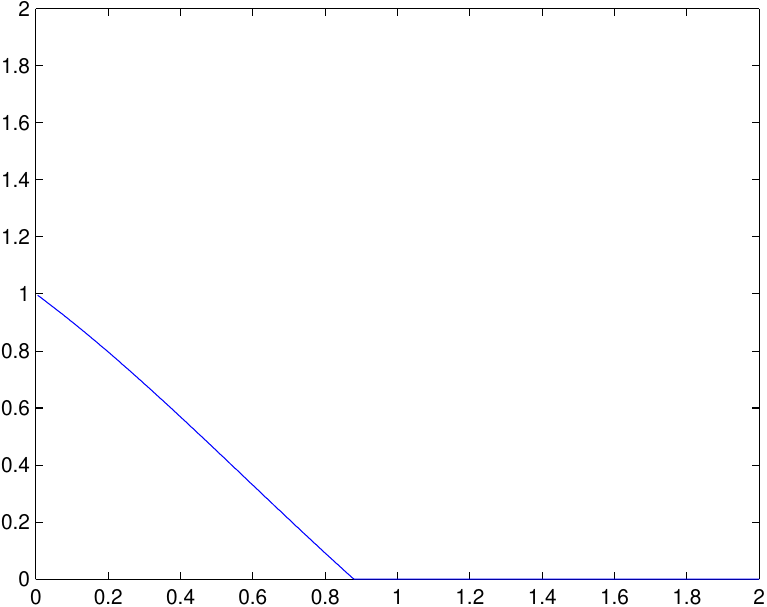}

\caption{Upper bound on the gap radius, as $ \Vert [V,U] \Vert $ varies from $0$ to about $0.21$.}
\label{fig:Gap-Upper-bound}
\end{figure}

An easy upper bound on the gap radius can be found, using the example in Lemma~\ref{lem:CyclicShift}, shown in
Fig.~\ref{fig:Gap-Upper-bound}.
This shows we cannot get as big a~gap using $B(U,V)$ as was possible with $VUV^{*}U^{*}$, but that the situation is
likely not as bad as Fig.~\ref{fig:smallGap} indicates.

Fortunately $B(U,V)$ is readily computable since~$f$,~$g$ and~$h$ we chosen to have rather fast decay in their Fourier
coef\/f\/icients.
Thus we were able to replace~\eqref{eq:FourierSeriesFunctCalc} by the simpler evaluation of order-$5$ trig polynomials.
For applications to index studies, the following is the most useful.
We will later have a~version of this for the Pfaf\/f\/ian--Bott index.

\begin{Proposition}
Suppose $
 \Vert  [U,V ] \Vert \leq0.206007$.
If $ \Vert  [U_{1},V_{1} ] \Vert \leq0.206007$ and $ \kappa(U_{1},V_{1})\neq\kappa(U,V) $ then
\begin{gather*}
 \Vert U-U_{1} \Vert + \Vert V-V_{1} \Vert \geq \frac{1}{5}\sqrt{1-5 \Vert
 [U,V ] \Vert^{2}} + \frac{1}{5}\sqrt{1-5 \Vert  [U_{1},V_{1} ] \Vert^{2}}.
\end{gather*}
\end{Proposition}

Another consequence of Theorem~\ref{thm:measure_the_gap} is the following.

\begin{Theorem}
Suppose~$U$ and~$V$ are unitary matrices.
If $ \Vert [U,V] \Vert \leq0.206007$ then
\begin{gather*}
\kappa(U,V)=\omega(U,V).
\end{gather*}
\end{Theorem}

\begin{proof}
Exel~\cite{ExelSoftTorusI} showed that for $ \Vert [U,V] \Vert <2$ the winding invariant equals an
abstract~$K$-theory invariant.
In the notation of~\cite{ExelSoftTorusI}, this is def\/ined in terms of $b_{\delta}$ in $K_{0}(A_{\delta})$.
It is easy to check that as long $B(u,v)$ has a~spectral gap, where~$u$ and~$v$ are the generators of the soft torus
$A_{\delta}$, the~$K$-theory class of $B(u,v)$ is $b_{\delta}$.
\end{proof}

\begin{Remark}
The smallest $\delta = \|[U,V]\|$ for which these invariants have been observed to dif\/fer, in numerical examples, is
$\delta \approx 0.85$.
\end{Remark}

\begin{Remark}
There is ambiguity in how we def\/ine the Bott index when $\|[U,V]\|$ is large.
Exel's abstract index in terms of the soft torus is not really computable, as there are no known explicit formulas for
$b_\delta$ when $\delta > 0.206007$.
Certainly the estimates above are not optimal so we are not sure for which~$\delta$ a~dif\/ferent formula is needed.
In general, when computing the~$K$-theoretical obstructions to approximate representations of relations being close to
exact representations, there is ambiguity when the error in the relations is large.
Ultimately, the ambiguity does not generally matter, as was discussed at length
in~\cite{EilersLoringContingenciesStableRelations}.
In the case of almost commuting unitary matrices, both the winding number invariant and the Bott index as in
Def\/inition~\ref{defn:bott_index} give computable invariants that, when the commutator is small, are stable in a~sizable
region around the pair.
More importantly, when either index is non-trivial, we know there is a~considerable distance to any pair of commuting
unitary matrices.
\end{Remark}

\section{Def\/ining the Pfaf\/f\/ian--Bott index}\label{sec:PfaffBott}

The Pfaf\/f\/ian of skew-symmetric matrices is not the most familiar object, and it it not clear at the outset how it
applies to a~problem involving self-dual matrices.
Let us start by recalling the dual operation.

We f\/ix
\begin{gather*}
Z=Z_{N}=\left[
\begin{matrix} 0 & I
\\
-I & 0
\end{matrix}
\right]
\end{gather*}
in $\mathbf{M}_{2N}(\mathbb{C})$ and this specif\/ies the dual operation
\begin{gather*}
X^{\sharp}=-ZX^{\mathrm{T}}Z
\end{gather*}
as above in~\eqref{eq:dualDefn}.

When we discuss $\mathbf{M}_{2}(\mathbf{M}_{2N}(\mathbb{C})) =
\mathbf{M}_{2N}(\mathbb{C})\otimes\mathbf{M}_{2}(\mathbb{C})$ we require the unitary
\begin{gather*}
Q=\frac{1}{\sqrt{2}}\left[
\begin{matrix} I & -iZ
\\
iZ & I
\end{matrix}
\right]=\frac{1}{\sqrt{2}}\left[
\begin{matrix} I & 0 & 0 & -iI
\\
0 & I & iI & 0
\\
0 & iI & I & 0
\\
-iI & 0 & 0 & I
\end{matrix}
\right]
\end{gather*}
which has the convenient property~\cite[Lemma 1.3]{HastLorTheoryPractice}
\begin{gather*}
Q^{*}X^{\sharp\otimes\sharp}Q= (Q^{*}XQ )^{\mathrm{T}}.
\end{gather*}
Here
\begin{gather*}
\left[
\begin{matrix} A & B
\\
C & D
\end{matrix}
\right]^{\sharp\otimes\sharp}=\left[
\begin{matrix} D^{\sharp} & -B^{\sharp}
\\
-C^{\sharp} & A^{\sharp}
\end{matrix}
\right]
\end{gather*}
or
\begin{gather*}
\left[
\begin{matrix} A_{11} & A_{12} & A_{13} & A_{14}
\\
A_{21} & A_{22} & A_{23} & A_{24}
\\
A_{31} & A_{32} & A_{33} & A_{34}
\\
A_{41} & A_{42} & A_{43} & A_{44}
\end{matrix}
\right]^{\sharp\otimes\sharp}=\left[
\begin{matrix} A_{44}^{\mathrm{T}} & -A_{34}^{\mathrm{T}} & -A_{24}^{\mathrm{T}} & A_{14}^{\mathrm{T}}
\vspace{1mm}\\
-A_{43}^{\mathrm{T}} & A_{33}^{\mathrm{T}} & A_{23}^{\mathrm{T}} & -A_{13}^{\mathrm{T}}
\vspace{1mm}\\
-A_{42}^{\mathrm{T}} & A_{32}^{\mathrm{T}} & A_{22}^{\mathrm{T}} & -A_{12}^{\mathrm{T}}
\vspace{1mm}\\
A_{41}^{\mathrm{T}} & -A_{31}^{\mathrm{T}} & -A_{21}^{\mathrm{T}} & A_{11}^{\mathrm{T}}
\end{matrix}
\right].
\end{gather*}

Recall the Pfaf\/f\/ian is def\/ined for all skew-symmetric, complex $2n$-by-$2n$ matrices~by
\begin{gather*}
\mathrm{Pf}\left(O\left[
\begin{matrix} 0 & a_{1}
\\
-a_{1} & 0 & a_{2}
\\
& -a_{2} & 0 & a_{3}
\\
& & -a_{3} & 0 & \ddots
\\
& & & \ddots & \ddots
\end{matrix}
\right]O^{\mathrm{T}}\right)=\det(O)a_{1}a_{3}\cdots a_{2n-1}
\end{gather*}
for~$O$ real orthogonal.
(All skew-symmetric matrices have such a~factorization, a~modif\/ied Hessenberg decomposition.) The essential properties
are that
\begin{gather*}
\mathrm{Pf}\big(YXY^{\mathrm{T}}\big)=\det(Y)\mathrm{Pf}(X)
\end{gather*}
for arbitrary~$Y$, that the Pfaf\/f\/ian varies continuously, and{\samepage
\begin{gather*}
 (\mathrm{Pf}(X) )^{2}=\det(X)
\end{gather*}
so the Pfaf\/f\/ian is zero exactly on the set of skew-symmetric, singular matrices.}

For matrices with the symmetry $X^{\sharp\otimes\sharp}=-X$ we can def\/ine a~modif\/ied Pfaf\/f\/ian
\begin{gather*}
\widetilde{\mathrm{Pf}}(X)=\mathrm{Pf} (Q^{*}XQ ).
\end{gather*}
We still have
\begin{gather*}
\big(\widetilde{\mathrm{Pf}}(X)\big)^{2}=\det(X)
\end{gather*}
and that this varies continuously.
The sign of the Pfaf\/f\/ian can be used to prove a~homotopy result, in the same way we use the determinant to detect that
the real orthogonal matrices fall into two connected parts.

\begin{Proposition}
Suppose~$B$ is in $\mathbf{M}_{4N}(\mathbb{C})$ and $B^{*}=B$ and $B^{\mathrm{T}}=-B$ and~$B$ is invertible.
Then $ \mathrm{Pf}(B)\in\mathbb{R}\setminus\{0\}.
$ If $B_{1}$ and $B_{2}$ are elements of
\begin{gather*}
\mathcal{H}=\big\{ B\in\mathbf{M}_{4N}(\mathbb{C})\,\big|\,B^{*}=B=-B^{\mathrm{T}}\; \text{is invertible}\big\},
\end{gather*}
then they can be connected by a~path in $\mathcal{H}$ if and only if with $\mathrm{Pf}(B_{1})$ and $\mathrm{Pf}(B_{2})$
have the same sign.
\end{Proposition}

\begin{proof}
We can apply Theorem~8.7 in~\cite{HastLorTheoryPractice} to $iB$ and learn that there is a~real orthogonal matrix~$O$ so
that $\det(O)=1$ and
\begin{gather*}
B=O\left[
\begin{matrix} 0 & i\lambda_{1}
\\
-i\lambda_{1} & 0
\\
& & 0 & i\lambda_{2}
\\
& & -i\lambda_{2} & 0
\\
& & & & \ddots & \ddots
\\
& & & & \ddots & \ddots
\end{matrix}
\right]O^{\mathrm{T}}
\end{gather*}
and all the real numbers $\lambda_{j}$ are positive except $\lambda_{1}$ which as the same sign as $\mathrm{Pf}(B)$.
It is clear from this form that two such matrices with Pfaf\/f\/ian of the same sign will be connected.
The spectrum of~$B$ is $ \{\pm\lambda_{1},\dots,\pm\lambda_{2N} \} $.
Thus there will be an even number of negative eigenvalues, so $\det(B)>0$.
Since the square of the Pfaf\/f\/ian is the determinant, we f\/ind $\mathrm{Pf}(B)$ is real, and as~$B$ is invertible, the
Pfaf\/f\/ian cannot be zero.
Since the Pfaf\/f\/ian varies continuously, it is not possible to connect two matrices in $\mathcal{H}$ that have Pfaf\/f\/ians
of opposite signs.
\end{proof}

Now we explain the Pfaf\/f\/ian--Bott index.

\begin{Definition}
Let~$f$,~$g$,~$h$ and $B(U,V)$ be as in Section~\ref{sec:Bott-Index}.
If $ \Vert [U,V] \Vert \leq0.206007$ def\/ine $\kappa_2(U,V)$ as the value in $\{\pm1\}$ given~by
\begin{gather*}
\kappa_{2}(U,V) = \mathrm{Sign}\big(\widetilde{\mathrm{Pf}} (B(U,V) )\big).
\end{gather*}
\end{Definition}

\begin{Lemma}
When~$U$ and~$V$ are commuting unitary matrices, $\kappa_2(U,V)=1$.
\end{Lemma}

\begin{proof}
This is a~special case of Theorem~8.8 of~\cite{HastLorTheoryPractice}.
Here is that proof in less technical language for this special case.

One easily checks that $\kappa_{2}$ remains constant along a~path so long as $ \Vert
[U_{t},V_{t}] \Vert \leq0.206007$.
One can use joint functional calculus (for commuting normal operators), and so keep the self-dual condition, to deform
a~commuting pair $U_{0}$ and $V_{0}$ that is self-dual over to $U_{1}=I$ and $V_{1}=I$.
One can then compute that
\begin{gather*}
B(I,I)=\left(
\begin{matrix} 0 & I
\\
I & 0
\end{matrix}
\right)
\end{gather*}
and
\begin{gather*}
\widetilde{\mathrm{Pf}}\left(
\begin{matrix} 0 & I
\\
 I & 0
\end{matrix}
\right)=\mathrm{Pf}\left[
\begin{matrix} iZ & 0
\\
0 & -iZ
\end{matrix}
\right] = \mathrm{Pf}\left(iZ\right)\mathrm{Pf}\left(-iZ\right) = i^{n}(-i)^{n} = 1.\tag*{\qed}
\end{gather*}
  \renewcommand{\qed}{}
\end{proof}

\begin{Proposition}
Suppose $
\left\Vert \left[U,V\right]\right\Vert \leq0.206007$ and that~$U$,~$V$, $U_{1}$, $V_{1}$ are self-dual unitary matrices.
If $\left\Vert \left[U_{1},V_{1}\right]\right\Vert \leq0.206007$ and $ \kappa_{2}(U_{1},V_{1})\neq\kappa_{2}(U,V) $ then
\begin{gather*}
\left\Vert U-U_{1}\right\Vert +\left\Vert V-V_{1}\right\Vert \geq \frac{1}{5}\sqrt{1-5\left\Vert
\left[U,V\right]\right\Vert^{2}} + \frac{1}{5}\sqrt{1-5\left\Vert \left[U_{1},V_{1}\right]\right\Vert^{2}}.
\end{gather*}

\end{Proposition}

\begin{Proposition}
Suppose $\Vert [U,V]\Vert \leq 0.206007$ and that~$U$, $V$, $U_{1}$, $V_{1}$ are self-dual unitary matrices.
If $U_{1}$ commutes with $V_{1}$ and $\kappa_{2}(U,V)=-1$ then
\begin{gather*}
\Vert U-U_{1} \Vert + \Vert V-V_{1} \Vert \geq\frac{1}{5} + \frac{1}{5}\sqrt{1-5 \Vert
[U,V]\Vert^{2}}.
\end{gather*}
\end{Proposition}

\begin{Remark}
We can describe the Pfaf\/f\/ian--Bott index more generally if we use the language of real $C^{*}$-algebras.
Suppose we are given a~real $C^{*}$-algebra as a~complex $C^{*}$-algebra~$B$ along with an anti-multiplicative linear
involution~$\tau$ on~$B$.
Then a~nice picture of $K_{-2}(B,\tau)$ (ignoring details with higher matrices) is in terms of
\begin{gather*}
\big\{x\in\mathrm{GL} (B\otimes\mathbf{M}_{2}(\mathbb{C}) ) \, \big| \, x^{*}=x,\; x^{\tau\otimes\sharp}=-x
\big\},
\end{gather*}
as was proven in~\cite[\S~8]{HastLorTheoryPractice}.
Then $B(z,w)$ is an element in $C(\mathbb{T}^{2})\otimes\mathbf{M}_{2}(\mathbb{C})$ that is hermitian, with spectrum
$\{-1,1\}$ and, with~$\tau$ the identity on $C(\mathbb{T}^{2})$, also $ B(z,w)^{\tau\otimes\sharp}=-B(z,w) $.
Therefore $B(z,w)$ determines an element in $K_{-2}(C(\mathbb{T}^{2}))$.
Given~$U$ and~$V$ self-dual unitary matrices in some unital real $C^{*}$-algebra $(B,\tau)$, we can def\/ine, for
now informally, the push-forward by something that is ``almost a~morphism'' $\psi: (C(\mathbb{T}^{2}),\tau)
\to (B,\tau)$, to produce the element~$[B(U,V)] $ in~$ K_{-2}(B,\tau).
$ We can def\/ine a~real structure on the soft-torus $A_{\delta}$ by setting $u_{\delta}^{\tau}=u_{\delta}$ and
$v_{\delta}^{\tau}=v_{\delta}$ (see~\cite[Chapter~5]{sorensen2012SPandGraphs} for details on why this is well-def\/ined) and consider
$\beta_{\delta}=[B(u_{\delta},v_{\delta})]$.
Then, for small~$\delta$, we have a~real $*$-homomorphism
\begin{gather*}
\gamma: \ (A_{\delta},\tau)\rightarrow (B,\tau)
\end{gather*}
and the Pfaf\/f\/ian--Bott element is then $\gamma_{*}(\beta_{\delta})$ in $K_{-2}(B,\tau)$.
For larger~$\delta$ we can proceed, but need an analysis of $K_{-2}(A_{\delta},\tau)$ that is best left for another
paper.
In the special case of $(B,\tau)$ equal to $(\mathbf{M}_{2n}(\mathbb{C}),\sharp)$ we used f\/irst
the isomorphism
\begin{gather*}
K_{-2} (\mathbf{M}_{2n}(\mathbb{C}),\sharp )\cong K_{2} (\mathbf{M}_{2n+2}(\mathbb{C}),\mathrm{T} )
\end{gather*}
induced by conjugation by a~set unitary, and then the isomorphism
\begin{gather*}
K_{2} (\mathbf{M}_{2n+2}(\mathbb{C}),\mathrm{T} )\rightarrow\mathbb{Z}/2
\end{gather*}
induced by the sign of the Pfaf\/f\/ian.
So
\begin{gather*}
\left(C(\mathbb{T}^{2}),\mathrm{id}\right)\leftarrow\left(A_{\delta},\tau\right)\rightarrow(\mathbf{M}_{2n}(\mathbb{C}),\sharp)
\end{gather*}
leads to
\begin{gather*}
K_{-2} (C(\mathbb{T}^{2}),\mathrm{id} )\leftarrow K_{-2} (A_{\delta},\tau )\rightarrow
K_{-2}(\mathbf{M}_{2n}(\mathbb{C}),\sharp)\rightarrow
K_{2} (\mathbf{M}_{2n+2}(\mathbb{C}),\mathrm{T} )\rightarrow\mathbb{Z}/2
\end{gather*}
with the right-most arrow given by the Pfaf\/f\/ian.
The remaining issue is showing the left-most arrow is surjective, which we have done here for small~$\delta$~by
explicitly def\/ining $B(u_\delta,v_\delta)$.
\end{Remark}

\section{Proof that the gap persists}
\label{sec:proof_of_gap}

Now we prove Theorem~\ref{thm:measure_the_gap}, f\/inding a~lower bound on the size of the gap in $B(U,V)$ as long as
$\delta = \| [U,V] \|$ is not too big.
We do so by f\/inding an upper bound on the norm of $B(U,V)^{2}-I$.
It is then a~routine application of the spectral mapping theorem to get lower bound on the size of the gap.

We will need some results about commutators and the functional calculus.
There is the folklore estimate $\Vert [f[V],U] \Vert \leq  \Vert
f^{\prime} \Vert_{\mathrm{F}} \Vert [U,V]\Vert $ where $\Vert
f^{\prime}\Vert_{\mathrm{F}} $ is the $\ell^{1}$-norm of the sequence of Fourier coef\/f\/icients of~$f$.
On its own, this estimate is really only helpful for very small commutators.

\begin{Definition}
Suppose~$f$ is continuous and $2\pi$-periodic.
Following~\cite{LoringVidesCommutators} we def\/ine $\eta_{f}:[0,\infty)\rightarrow[0,\infty)$~by
\begin{gather*}
\eta_{f}(\delta) = \sup \big\{ \Vert  [f[V],A ] \Vert \,| \, V~\mbox{is unitary},\; \Vert A\Vert \leq1, \; \Vert [V,A]\Vert \leq\delta \big\},
\end{gather*}
where the supremum is taken over all~$V$ and~$A$ in every unital $C^{*}$-algebra.

Once we have a~bound on $\eta_{f}$ we can use it to bound more than just commutators.
Indeed, by~\cite[Lemma~1.2]{LoringVidesCommutators}, for any two unitaries $V_{1}$and $V_{2}$ we have
\begin{gather*}
 \Vert f[V]-f[V_{1}] \Vert \leq\eta_{f} ( \Vert V-V_{1} \Vert ).
\end{gather*}
We need a~special case of a~lemma in~\cite{LoringVidesCommutators}.
\end{Definition}

\begin{Lemma}
\label{lem:eta_Bounds_m_delta_plus_b}
Suppose~$f$ is continuous, real-valued and periodic, and that $f_1$ is the trigonometric polynomial
\begin{gather*}
f_1(x)=\sum\limits_{k=-n}^{n}a_{k}e^{ikx}.
\end{gather*}
Let $f_2 = f - f_1$.
Then $\eta_f(\delta)\leq m\delta+b$ where
\begin{gather*}
m=\sum\limits_{k=-n}^{n} |ka_{k} |
\end{gather*}
and $ b=\max f_2(x)-\min f_2(x) $.
\end{Lemma}

Before we focus on our choice of the three functions $f$,~$g$ and~$h$ to use in the Bott invariant, we look at the terms
we need to control when bounding $B(U,V)^{2}-I$.

\begin{Lemma}
Suppose $f$,~$g$ and~$h$ are continuous, real-valued functions that are $2\pi$-periodic, and with $ f^{2}+g^{2}+h^{2}=1
$ and $ gh=0$. Suppose~$U$ and~$V$ are unitary matrices and define
\begin{gather*}
S=\left[
\begin{matrix} f[V] & g[V]+\tfrac{1}{2} \{h[V],U \}
\vspace{1mm}\\
g[V]+\tfrac{1}{2} \{h[V],U^{*} \} & -f[V]
\end{matrix}
\right].
\end{gather*}
Then $S^{*}=S$ and
\begin{gather*}
 \big\Vert S^{2}-I \big\Vert \leq2 \Vert  [h[V],U ] \Vert + \Vert  [f[V],U ] \Vert.
\end{gather*}
\end{Lemma}

\begin{proof}
Since~$f$,~$g$ and~$h$ are real-valued, the matrices $f[V]$, $g[V]$ and $h[V]$ are hermitian.
Let us write~$f$ for $f[V]$, etc.
We see easily $S^{*}=S$ and
\begin{gather*}
\left[
\begin{matrix} f & g+\tfrac{1}{2} \{h,U \}
\\
g+\tfrac{1}{2}\left\{h,U^{*}\right\} & -f
\end{matrix}
\right]^{2}-\left[
\begin{matrix} I & 0
\\
0 & I
\end{matrix}
\right]=\left[
\begin{matrix} A & B
\\
B^{*} & A^{*}
\end{matrix}
\right],
\end{gather*}
where
\begin{gather*}
A   =f^{2}+g^{2}-I+\tfrac{1}{4} \{h,U \}  \{h,U^{*} \} +\tfrac{1}{2}g \{h,U^{*} \}
+\tfrac{1}{2} \{h,U \} g
\\
\hphantom{A}
  =-h^{2}+\tfrac{1}{4} \{h,U \}  \{h,U^{*} \} +\tfrac{1}{2}g \{h,U^{*} \}
+\tfrac{1}{2} \{h,U \} g
\end{gather*}
and
\begin{gather*}
B =fg+\tfrac{1}{2}f \{h,U \} -gf-\tfrac{1}{2} \{h,U \} f =\tfrac{1}{2}f \{h,U \}
-\tfrac{1}{2} \{h,U \} f.
\end{gather*}
We have
\begin{gather*}
\big\Vert S^{2}-I\big\Vert   \leq\left\Vert \left[
\begin{matrix} A & 0
\\
0 & A^{*}
\end{matrix}
\right]\right\Vert +\left\Vert \left[
\begin{matrix} 0 & B
\\
B^{*} & 0
\end{matrix}
\right]\right\Vert
  =\Vert A\Vert +\Vert B\Vert.
\end{gather*}
Notice $f^{2}+g^{2}+h^{2}=1$ forces these functions to take value in $[-1,1]$ so $\Vert f[V]\Vert \leq1$,
etc.
Therefore
\begin{gather*}
 \Vert A \Vert   \leq\tfrac{1}{4}\big\Vert hUhU^{*}+Uh^{2}U^{*}+UhU^{*}h-3h^{2}\big\Vert
+\tfrac{1}{2} \Vert gU^{*}h-ghU^{*}+hUg-Uhg \Vert
\\
\hphantom{\Vert A \Vert}{}
  \leq\tfrac{1}{2} \Vert h \Vert  \Vert Uh-hU \Vert +\tfrac{1}{4}\big\Vert Uh^{2}-h^{2}U\big\Vert
+ \Vert g \Vert  \Vert Uh-hU \Vert
  \leq2 \Vert Uh-hU \Vert
\end{gather*}
and
\begin{gather*}
 \Vert B \Vert   =\tfrac{1}{2} \Vert hfU+fUh-hUf-Ufh \Vert
  \leq\tfrac{1}{2} \Vert h [f,U ] \Vert +\tfrac{1}{2} \Vert  [f,U ]h \Vert
  \leq \Vert [f,U ] \Vert,
\end{gather*}
so
\begin{gather*}
\big\Vert S^{2}-I\big\Vert \leq2 \Vert  [h,U ] \Vert + \Vert  [f,U ] \Vert.\tag*{\qed}
\end{gather*}
  \renewcommand{\qed}{}
\end{proof}

Now we let~$f$,~$g$ and~$h$ be the functions from Def\/inition~\ref{Def:BoffByf_g_h}.
Here we start needing a~computer algebra package.
It shows us that
\begin{gather*}
f(x)^{2}+\frac{407}{512}\cos^{6}(x)\left(1+\frac{96}{407}\cos (2x )+\frac{9}{407}\cos (4x )\right)=1,
\end{gather*}
which means
\begin{gather*}
g(x)=\sqrt{\frac{407}{512}}\cos^{3}(x)\sqrt{1+\frac{96}{407}\cos(2x)+\frac{9}{407}\cos(4x)}\left(1-\chi_{[-\frac{\pi}{2},\frac{\pi}{2}]}(x)\right)
\end{gather*}
and
\begin{gather*}
h(x)=\sqrt{\frac{407}{512}}\cos^{3}(x)\sqrt{1+\frac{96}{407}\cos(2x)+\frac{9}{407}\cos(4x)}\chi_{[-\frac{\pi}{2},\frac{\pi}{2}]}(x).
\end{gather*}
A~handy formula here is
\begin{gather*}
\frac{407}{320}\left(1+\frac{96}{407}\cos(2x)+\frac{9}{407}\cos(4x)\right)=\left(1+\frac{15}{40}\cos^{2}(x)+\frac{9}{40}\cos^{4}(x)\right)
\end{gather*}
and we get alternate expression for~$g$ and~$h$, in particular
\begin{gather*}
h(x)=\frac{\sqrt{10}}{4}\cos^{3}(x)\sqrt{1+\frac{15}{40}\cos^{2}(x)+\frac{9}{40}\cos^{4}(x)}\chi_{[-\frac{\pi}{2},\frac{\pi}{2}]}(x).
\end{gather*}

We new bound the derivative of~$g$ and~$h$, computing
\begin{gather*}
\frac{d}{dx}\left(\frac{\sqrt{10}}{4}\cos^{3}(x)\sqrt{1+\frac{15}{40}\cos^{2}(x)+\frac{9}{40}\cos^{4}(x)}\right)
  =\frac{p(\sin(x))}{16\sqrt{q(\sin(x))}},
\end{gather*}
where
\begin{gather*}
p(x)=30x(x-1)(x+1)\big(3x^{4}-10x^{2}+15\big)
\end{gather*}
and
\begin{gather*}
q(x)=9x^{4}-33x^{2}+64.
\end{gather*}
On $[-1,1]$ the max of $p(x)$ is $150$ and the min of $q(x)$ is $64$ so we f\/ind
\begin{gather*}
|h^{\prime}(x)|\leq\frac{150}{128}
\end{gather*}
and the same for $g^{\prime}$.

We need~$h$ as a~Fourier series so need
\begin{gather*}
c_{n}=\frac{1}{2\pi}
\int_{-\frac{\pi}{2}}^{\frac{\pi}{2}}\cos(nx)\sqrt{\frac{407}{512}}\cos^{3}(x)\sqrt{1+\frac{96}{407}\cos(2x)+\frac{9}{407}\cos(4x)}\,
dx.
\end{gather*}
We computed these with numerical integration, and without checking error estimates, in~\cite{LorHastHgTe}.
We compute these a~little more carefully here.
Thus Table~\ref{tab:coeff_of_f_g_h} is a~slightly more accurate replacement for Table~11
in~\cite{HastLorTheoryPractice}.
We f\/ind
\begin{gather*}
c_{n}
=\frac{1}{2\pi}\int_{-\frac{\pi}{2}}^{\frac{\pi}{2}}\cos(nx)\sqrt{\frac{407}{512}}\cos^{3}(x)\sum\limits_{k=0}^{\infty}{1/2
\choose k}\left(\frac{96}{407}\cos(2x)+\frac{9}{407}\cos(4x)\right)^{k}  dx
\\
\hphantom{c_{n}}{}
  =\sum\limits_{k=0}^{\infty}\frac{1}{2\pi}\sqrt{\frac{407}{512}}{1/2 \choose
k}\! \int_{-\frac{\pi}{2}}^{\frac{\pi}{2}}\cos(nx)\cos^{3}(x)\left(\frac{96}{407}\cos(2x)+\frac{9}{407}\cos(4x)\right)^{k}
dx
  =\sum\limits_{k=0}^{\infty}I_{n,k},
\end{gather*}
where the $I_{n,k}$ were def\/ined in-line and are easy to compute with a~computer algebra package.
The convergence here is rather rapid, as
\begin{gather*}
I_{n,k}   \leq\frac{1}{2\pi}\sqrt{\frac{407}{512}}{1/2 \choose
k}(-1)^{k}\int_{-\frac{\pi}{2}}^{\frac{\pi}{2}}\left|\cos(nx)\cos^{3}(x)\left(\frac{96}{407}\cos(2x)+\frac{9}{407}\cos(4x)\right)^{k}\right|
dx
\\
\hphantom{I_{n,k}}{}
  \leq\frac{1}{2\pi}\sqrt{\frac{407}{512}}{1/2 \choose k}\left(-\frac{105}{407}\right)^{k}.
\end{gather*}
Letting $T_{K}$ denote the Taylor polynomial $ T_{K}(x)\approx\sqrt{1+x} $ of degree~$K$ expanded at~$0$, we have
\begin{gather*}
\sum\limits_{k=K+1}^{\infty} |I_{n,k} |
\leq\frac{1}{2\pi}\sqrt{\frac{407}{512}}\sum\limits_{k=K+1}^{\infty}{1/2 \choose k}\left(-\frac{105}{407}\right)^{k}
\\
\hphantom{\sum\limits_{k=K+1}^{\infty} |I_{n,k} |}{}
  =\frac{1}{2\pi}\sqrt{\frac{407}{512}}\left(\sum\limits_{k=0}^{\infty}{1/2 \choose
k}\left(-\frac{105}{407}\right)^{k}-\sum\limits_{k=0}^{K}{1/2 \choose k}\left(-\frac{105}{407}\right)^{k}\right)
\\
\hphantom{\sum\limits_{k=K+1}^{\infty} |I_{n,k} |}{}
  =\frac{1}{2\pi}\sqrt{\frac{407}{512}}\left(\sqrt{1-\frac{105}{407}}-T_{K}\left(-\frac{105}{407}\right)\right).
\end{gather*}
This means we need $K=7$ to get six digits absolute accuracy, with the results shown in Table~\ref{tab:coeff_of_f_g_h}.
The integration was done symbolically in Matlab\footnote{Code assisting with tables and f\/igures and calculations is
available at \url{http://repository.unm.edu/handle/1928/23494}.}.

\begin{table}[t]\centering
\caption{These are approximations to the f\/irst coef\/f\/icients in the Fourier expansions of the~$f$,~$g$ and~$h$ used to
def\/ine the Bott index.
Extend these to negative indices by the rules $a_{-n}=\overline{a_{n}}$ and $b_{-n}=b_{n}$ and $c_{n}=-c_{n}$.}\label{tab:coeff_of_f_g_h}
\vspace{1mm}

\begin{tabular}{|c|c|c|c|c|c|c|}
\hline
$n$ & $0$ & $1$ & $2$ & $3$ & $4$ & $5$\\
\hline
\hline
$a_{n}$ & 0 & $-\frac{{150i}}{256}$ & 0 & $-\frac{{25i}}{256}$ & 0 &
$-\frac{{3i}}{256}$\tsep{2pt}\bsep{2pt}\\
\hline
$b_{n}$ & $0.202047$ & $-0.179940$ & $0.125655$ & $-0.066010$ & $0.023445$ & $-0.003886$\\
\hline
$c_{n}$ & $0.202047$ & $0.179940$ & $0.125655$ & $0.066010$ & $0.023445$ & $0.003886$\\
\hline
\end{tabular}
\end{table}

Using the values in the table to def\/ine
\begin{gather*}
h_{5}(x)=\sum\limits_{n=-5}^{5}c_{n}e^{inx}
\end{gather*}
we f\/ind
\begin{gather*}
\left|\frac{d}{dx} (h-h_{5} )\right|\leq\frac{150}{128}+1.48498=2.656855.
\end{gather*}
and so we can estimate to six decimal places the maximum of $ |h-h_{5} |$ by simply plugging in values between
$-\pi$ and~$\pi$ with an even spacing of a~little less than~$10^{-7}$.
Keeping track of the errors and rounding up, we f\/ind
\begin{gather*}
\mathrm{diam} (h(x)-h_{5}(x) )\leq0.004110
\end{gather*}
and we note
\begin{gather*}
\left\Vert h_{5}^{\prime}\right\Vert_{\rm F}=\sum\limits_{n=-5}^{5}n|c_{n}|=1.48498.
\end{gather*}

The other estimates of this sort, for $h_{0},\dots,h_{4}$, are summarized in Table~\ref{tab:forEta_h}.
We also can use brute force to f\/ind
\begin{gather*}
\sup_{x} |h(x)-h_{5}(x) |\leq0.002338.
\end{gather*}

\begin{Lemma}
For any unitary matrix~$V$,
\begin{gather*}
 \Vert h_{5}[V]-h[V] \Vert \leq0.002338.
\end{gather*}
\end{Lemma}

We get the same error estimate on using only $b_{-5}$ through $b_{5}$ when numerically compu\-ting~$g[V]$.

\begin{table}[t]\centering
\caption{Bounds on $\eta_{h}$ as a~slope and an of\/fset.}\label{tab:forEta_h}
\vspace{1mm}
\begin{tabular}{|c|c|c|}
\hline
$n$ & $m=$ bound on $ \Vert h_{n}^{\prime} \Vert_{\rm F}$ & $b=$ bound on
$\mathrm{diam}(h(x)-h_{n}(x))$\bsep{1pt}\tsep{1pt}\\
\hline
\hline
$0$ & $0$ & $1$\\
\hline
$1$ & $0.359880$ & $0.732237$\\
\hline
$2$ & $0.862500$ & $0.350141$\\
\hline
$3$ & $1.258560$ & $0.106619$\\
\hline
$4$ & $1.446120$ & $0.017509$\\
\hline
$5$ & $1.48498$ & $0.004110$\\
\hline
 $\infty$ & $2.99208$ & $0$\\
\hline
\end{tabular}
\end{table}

\begin{Lemma}
For~$h$ as in Definition~{\rm \ref{Def:BoffByf_g_h}}, we have
\begin{gather*}
 \Vert h^{\prime} \Vert_{\mathrm{F}}\leq2.99208.
\end{gather*}
\end{Lemma}

\begin{proof}
We check that
\begin{gather*}
h^{\prime}(x)=\frac{-\frac{1}{\sqrt{3256}}\big(45\cos^{4}(x)+60\cos^{2}(x)+120\big)}{\sqrt{1+\frac{96}{407}\cos(2x)
+\frac{9}{407}\cos(4x)}}\sin(x)\cos^{2}(x)\chi_{[-\frac{\pi}{2},\frac{\pi}{2}]}(x)
\end{gather*}
and attack this as three factors.
It is easy to see
\begin{gather*}
\left\Vert
-\frac{1}{\sqrt{3256}}\big(45\cos^{4}(x)+60\cos^{2}(x)+120\big)\right\Vert_{\mathrm{F}}=\frac{225}{\sqrt{3256}}
\end{gather*}
and the next factor is not so bad, as we see
\begin{gather*}
\left\Vert \frac{1}{\sqrt{1+\frac{96}{407}\cos(2x)+\frac{9}{407}\cos(4x)}}\right\Vert_{\rm F}
\leq\sum\limits_{k=0}^{\infty}\left|{-1/2 \choose k}\right|\left\Vert
\frac{96}{407}\cos(2x)+\frac{9}{407}\cos(4x)\right\Vert_{\rm F}^{k}
\\
\qquad{}  =\sum\limits_{k=0}^{\infty}{-1/2 \choose k}(-1)^{k}\left(\frac{105}{407}\right)^{k}
   =\frac{1}{\sqrt{1-\frac{105}{407}}}=\sqrt{\frac{407}{302}}.
\end{gather*}
We estimate
\begin{gather*}
\big\Vert \sin(x)\cos^{2}(x)\chi_{[-\frac{\pi}{2},\frac{\pi}{2}]}(x)\big\Vert_{\mathrm{F}}
\end{gather*}
as follows.
The Fourier series of $-i\sin(x)\cos^{2}(x)\chi_{[-\frac{\pi}{2},\frac{\pi}{2}]}$ is
\begin{gather*}
\dots,\frac{16}{3465\pi},0,\frac{-4}{315\pi},0,\frac{8}{105\pi},\frac{1}{16},\frac{-8}{15\pi},\frac{-1}{16},0,\frac{1}{16},\frac{8}{15\pi},\frac{-1}{16},\frac{-8}{105\pi},0,\frac{4}{315\pi},0,\frac{-16}{3465\pi},\dots
\end{gather*}
with terms beyond $n=3$ being given~by
\begin{gather*}
-\cos\left(\frac{\pi n}{2}\right)\left(\frac{1}{(n-1)^{3}-4(n-1)}+\frac{1}{(n+1)^{3}-4(n+1)}\right).
\end{gather*}
Therefore
\begin{gather*}
  \big\Vert \sin(x)\cos^{2}(x)\chi_{[-\frac{\pi}{2},\frac{\pi}{2}]}\big\Vert_{\rm F}
=\frac{1}{4}+\frac{16}{15\pi}+\frac{18}{105\pi}+\frac{4}{\pi}\sum\limits_{n=3}^{\infty}\frac{1}{\left(2n+1\right)^{3}-4\left(2n+1\right)}
\\
\qquad{}
\leq\frac{1}{4}+\frac{16}{15\pi}+\frac{18}{105\pi}+\frac{4}{315\pi}+\frac{16}{3465\pi}+\frac{4}{\pi}\int_{4}^{\infty}\frac{1}{8x^{3}+12x^{2}-2x-3}\,
dx
\\
\qquad{}
=\frac{1}{4}+\frac{16}{15\pi}+\frac{18}{105\pi}+\frac{4}{315\pi}+\frac{16}{3465\pi}+\frac{1}{4\pi}\ln\left(\frac{81}{77}\right)
\end{gather*}
and so
\begin{gather*}
\Vert
h^{\prime}\Vert_{\mathrm{F}}\leq\frac{225}{\sqrt{3256}}\sqrt{\frac{407}{302}}\!
\left(\frac{1}{4}+\frac{16}{15\pi}+\frac{18}{105\pi}+\frac{4}{315\pi}+\frac{16}{3465\pi}+\frac{1}{4\pi}\ln\left(\frac{81}{77}\right)\right)
< 2.992076.\!\!\!\!\!\!\tag*{\qed}
\end{gather*}
  \renewcommand{\qed}{}
\end{proof}

We approximate~$f$ the same way, but this is just arithmetic since~$f$ is already a~trigonometric polynomial.

\begin{table}[t]\centering
\caption{Bounds on $\eta_{f}$ as a~slope and an of\/fset.}\label{tab:forEta_f}
\vspace{1mm}

\begin{tabular}{|c|c|c|}
\hline
$n$ & $m=$ bound on $ \Vert f_{n}^{\prime} \Vert_{\rm F}$ & $b=$ bound on
$\mathrm{diam}(f(x)-f_{n}(x))$\tsep{1pt}\bsep{1pt}\\
\hline
\hline
$0$ & $0$ & $2$\\
\hline
$1$ & $1.171875$ & $0.4375$\\
\hline
$2$ & $1.7578125$ & $0.04687$\\
\hline
$\infty$ & $1.875$ & $0$\\
\hline
\end{tabular}

\end{table}

\begin{Lemma}
\label{lem:eta_f_eat_g_eta_h}
Let~$f$ and~$g$ and~$h$ be as in Definition~{\rm \ref{Def:BoffByf_g_h}}.
Then $ \eta_{f}(\delta)\leq m\delta+b $ for each of the values in Table~{\rm \ref{tab:forEta_f}}, and $ \eta_{g}(\delta)\leq
m\delta+b $ and $ \eta_{h}(\delta)\leq m\delta+b $ for each of the values in Table~{\rm \ref{tab:forEta_h}}.
\end{Lemma}

\begin{figure}[t]\centering
\includegraphics[clip,scale=0.5]{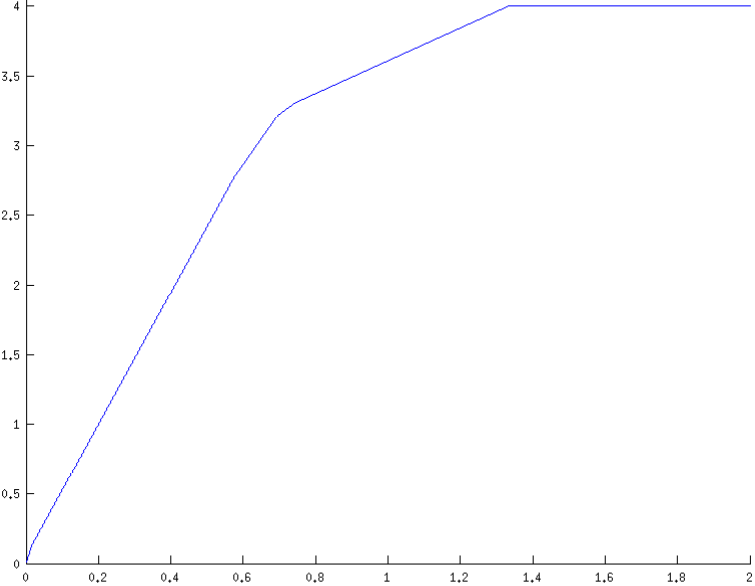}

\caption{The function $\beta(\delta)$ that bounds $ \Vert B(U,V)^{2}-I \Vert $ in terms of $\delta= \Vert
[U,V] \Vert $.}\label{fig:Bott_Bound_beta}
\end{figure}

Let $ \beta(\delta)=2\eta_{h}(\delta)+\eta_{f}(\delta) $ which is shown in Fig.~\ref{fig:Bott_Bound_beta}.

\begin{Theorem}
Suppose~$U$ and~$V$ are unitary matrices.
Then
\begin{gather*}
\big\Vert B(U,V)^{2}-I\big\Vert \leq\beta ( \Vert [U,V] \Vert  )
\end{gather*}
and for $ \Vert [U,V] \Vert \leq0.206007$ the gap at $0$ in the spectrum of $B(U,V)$ has radius at
least
\begin{gather*}
\sqrt{1-\beta ( \Vert [U,V] \Vert )}.
\end{gather*}
\end{Theorem}

The other key thing we must show is how $B(U,V)$ varies as~$U$ and~$V$ vary.
After this, all our main theorems will follow.

\begin{Theorem}
If $U_{j}$ and $V_{j}$ are unitary matrices then
\begin{gather*}
 \Vert B (U_{0},V_{0} )-B (U_1,V_1 ) \Vert \leq \beta ( \Vert V_{0}-V_{1} \Vert
 )+ \Vert U_{0}-U_{1} \Vert
\end{gather*}
and so
\begin{gather*}
 \Vert B (U_{0},V_{0} )-B (U_1,V_1 ) \Vert \leq \beta ( \Vert V_{0}-V_{1} \Vert
+ \Vert U_{0}-U_{1} \Vert  ).
\end{gather*}
\end{Theorem}

\begin{proof}
This follows easily from Lemma~\ref{lem:eta_f_eat_g_eta_h} and~\cite[Lemma~1.2]{LoringVidesCommutators}.
\end{proof}

\section{The log method}

An alternate way to compute the Bott index was considered in~\cite{ExelLoringInvariats}.
One replaces $B(U,V)$ with
\begin{gather*}
B_{\mathrm{L}}(U,V)=\left(
\begin{matrix} \frac{1}{\pi}K & \frac{1}{2}\left\{\sqrt{I-\frac{1}{\pi^{2}}K^{2}},U\right\}
\\
\frac{1}{2}\left\{\sqrt{I-\frac{1}{\pi^{2}}K^{2}},U^{*}\right\} & -\frac{1}{\pi}K
\end{matrix}
\right)
\end{gather*}
where $iK$ is the logarithm of~$V$, meaning $-\pi\leq K<\pi$ and $e^{iK}=V$.
Numerical evidence in~\cite{LorHastHgTe} suggests that, for small commutators, the Pfaf\/f\/ian--Bott index can be computed
using $B_{\mathrm{L}}(U,V)$.
We validate this here.

Since the logarithm is not continuous, numerical errors will mean we might accidentally compute the wrong branch of
logarithm on~$V$, or indeed any logarithm of~$V$ whatsoever.

We note that when~$q$ is periodic, $ q(K)=q[V] $.

\begin{figure}[t]\centering
\includegraphics[clip,scale=0.4]{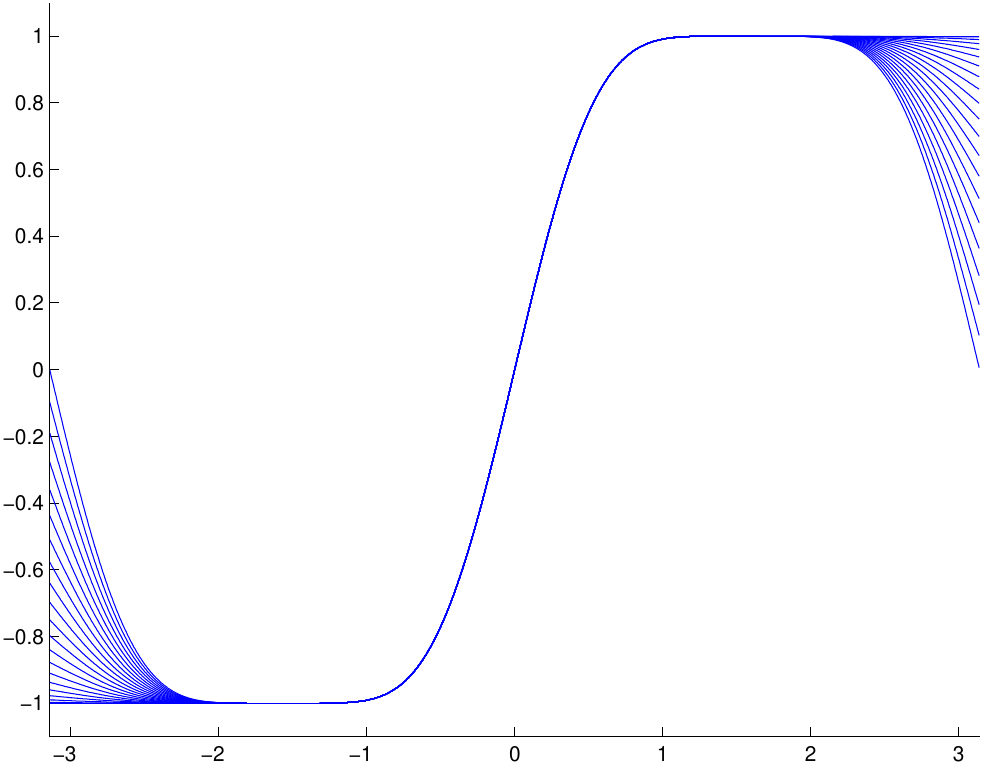}\qquad \includegraphics[clip,scale=0.4]{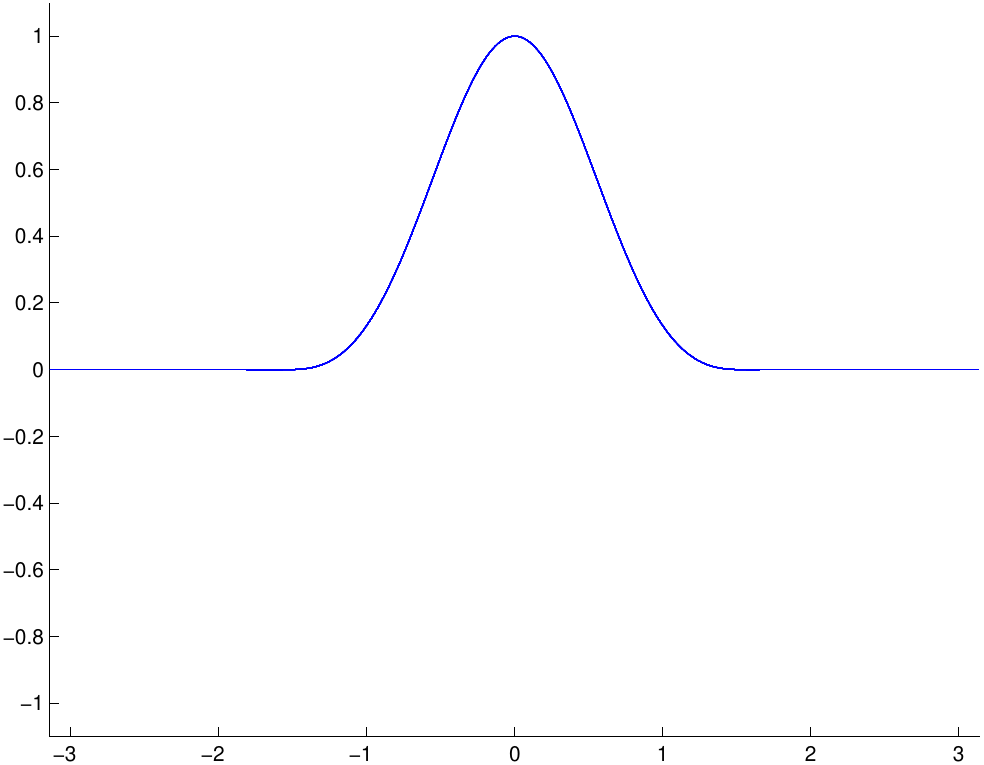}
\vspace{3mm}

\includegraphics[clip,scale=0.4]{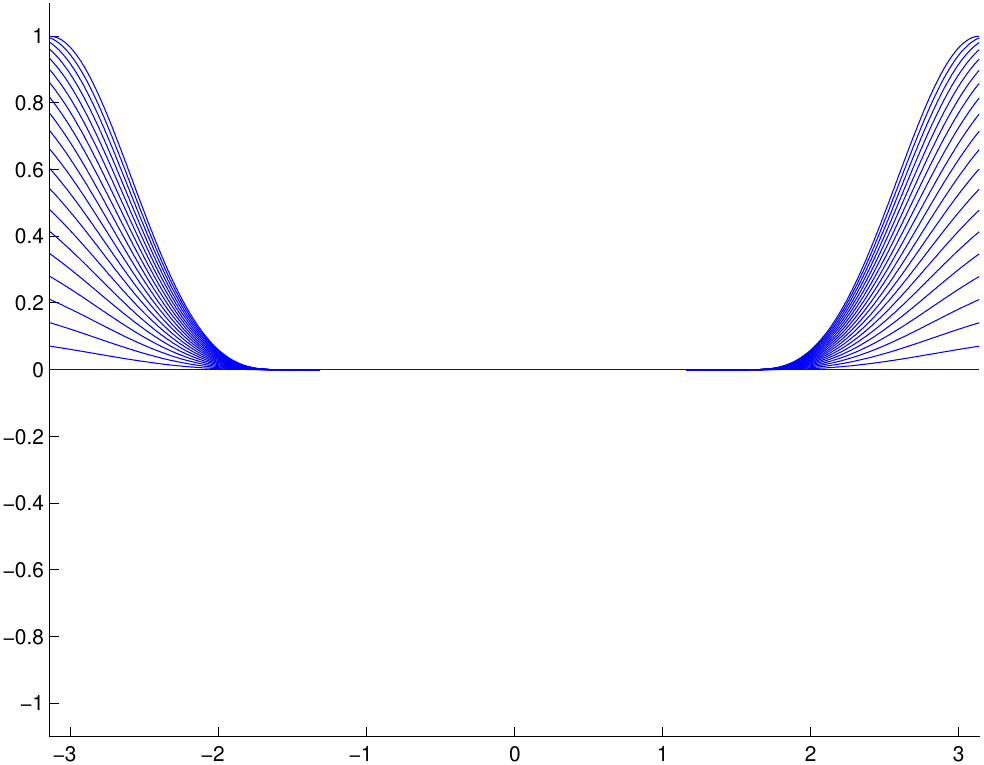}

\caption{First part of the path, showing $f_{t}$ then $h_{t}$ then $g_{t}$.}\label{Fig6}
\end{figure}

\begin{Lemma}\label{lem:borel_f_g_h}
Suppose $f$,~$g$ and~$h$ are real-valued Borel functions on $[-\pi,\pi]$ satisfying $ f^{2}+g^{2}+h^{2}=1 $ and $ gh=0
$.
Let $q(x)=f(x)h(x)$ and assume further that~$q$ and~$h$ are continuous and $2\pi$-periodic.
Suppose~$U$ and~$V$ are unitary matrices and $-iK$ is a~logarithm of~$V$ and define
\begin{gather*}
S=\left[
\begin{matrix} f(K) & g(K)+\tfrac{1}{2}\left\{h(K),U\right\}
\\
g(K)+\tfrac{1}{2}\left\{h(K),U^{*}\right\} & -f(K)
\end{matrix}
\right].
\end{gather*}
Then $S^{*}=S$ and
\begin{gather*}
\big\Vert S^{2}-I\big\Vert \leq  ( \Vert g \Vert +1 ) \Vert  [h[V],U ] \Vert +
\frac{1}{4} \Vert  [h[V],U ] \Vert^{2}+\frac{1}{2}\big\Vert \big[h^{2}[V],U\big]\big\Vert
+ \Vert  [q[V],U ] \Vert.
\end{gather*}
\end{Lemma}

\begin{proof}
We write~$f$ for $f(K)$, etc., and estimate a~bit more carefully than before.
We f\/ind
\begin{gather*}
  \frac{1}{4}\big\Vert hUhU^{*}+Uh^{2}U^{*}+UhU^{*}h-3h^{2}\big\Vert
\\
\qquad{}
  =\frac{1}{4}\big\Vert hUhU^{*}-h^{2}+UhU^{*}h-Uh^{2}U^{*}+2Uh^{2}U^{*}-2h^{2}\big\Vert
\\
  =\frac{1}{4}\big\Vert -[h,U][h,U]^{*}+2\big(Uh^{2}U^{*}-h^{2}\big)\big\Vert
\leq\frac{1}{4} \Vert  [h,U ] \Vert^{2}+\frac{1}{2}\big\Vert \big[h^{2},U\big]\big\Vert
\end{gather*}
and
\begin{gather*}
\frac{1}{2} \Vert g \{h,U^{*} \} + \{h,U \} g \Vert   =\frac{1}{2} \Vert
gU^{*}h+hUg \Vert
  =\frac{1}{2} \Vert g [U^{*},h ]+ [h,U ]g \Vert
  \leq \Vert g \Vert  \Vert  [h,U ] \Vert
\end{gather*}
and
\begin{gather*}
\tfrac{1}{2} \Vert f \{h,U \} - \{h,U \} f \Vert   =\tfrac{1}{2} \Vert
2fhU-2Ufh+fUh-fhU-hUf+Uhf \Vert
\\
\qquad{}  \leq \Vert fhU-Uhf \Vert +\tfrac{1}{2} \Vert fUh-fhU \Vert +\tfrac{1}{2} \Vert
hUf-Uhf \Vert
\\
\qquad {}
  \leq \Vert fhU-Uhf \Vert +\tfrac{1}{2} \Vert Uh-hU \Vert +\tfrac{1}{2} \Vert hU-Uh \Vert
\\
\qquad{}
  = \Vert  [fh,U ] \Vert + \Vert  [h,U ] \Vert
  = \Vert  [q,U ] \Vert + \Vert [h,U ] \Vert.\tag*{\qed}
\end{gather*}
  \renewcommand{\qed}{}
\end{proof}

\begin{figure}[t]\centering
\includegraphics[clip,scale=0.4]{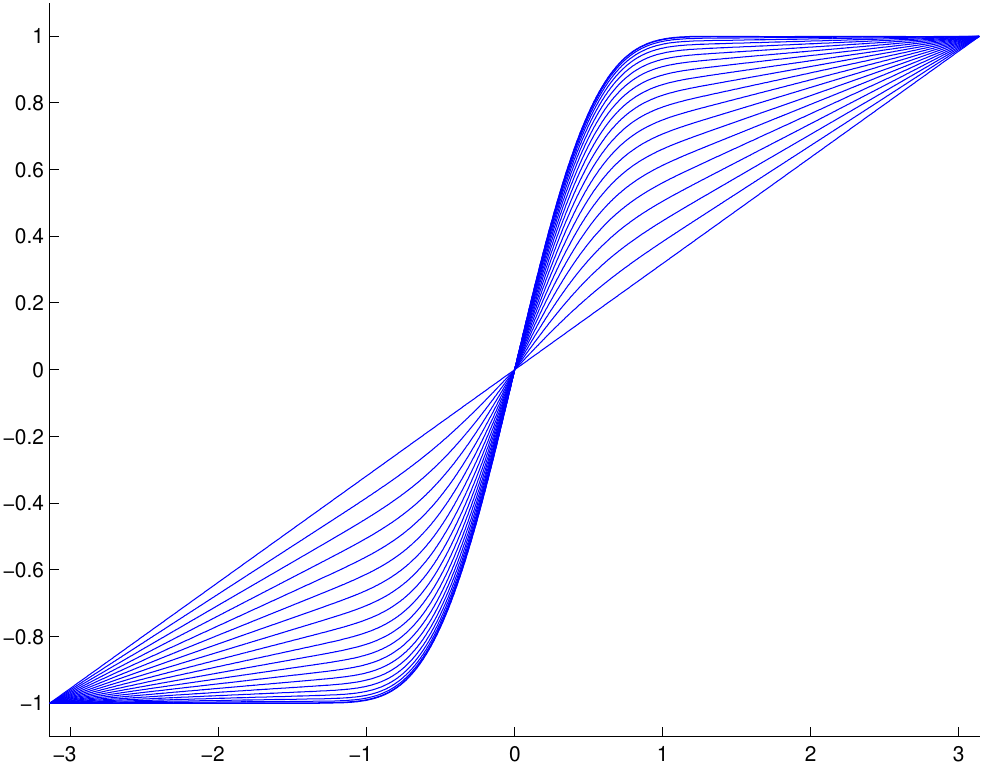}\qquad \includegraphics[clip,scale=0.4]{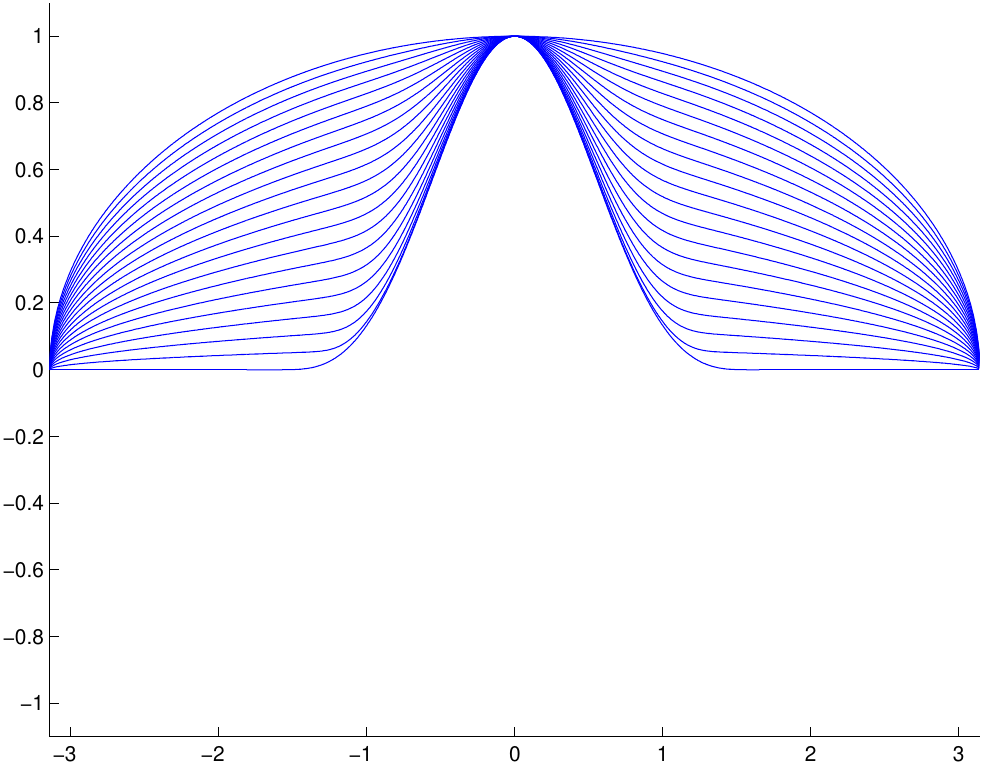}

\caption{Second
part of the path, showing only $f_{t}$ then $h_{t}$. Here $g_{t}$ is zero.}\label{Fig7}
\end{figure}

\begin{Lemma}
Suppose~$U$ and~$V$ are unitary matrices.
If $\left\Vert [U,V]\right\Vert \leq\frac{1}{8}$ then for any choice of~$K$ with $-\pi\leq K\leq\pi$ and
$e^{iK}=V$, there is a~path $B_{t}$ of invertible self-adjoint matrices between $B(U,V)$ and
\begin{gather*}
\left(
\begin{matrix} \frac{1}{\pi}K & \frac{1}{2}\left\{\sqrt{I-\frac{1}{\pi^{2}}K^{2}},U\right\}
\\
\frac{1}{2}\left\{\sqrt{I-\frac{1}{\pi^{2}}K^{2}},U^{*}\right\} & -\frac{1}{\pi}K
\end{matrix}
\right)
\end{gather*}
and, if~$U$ and~$V$ are self-dual, then the path may be chosen with the symmetry $B_{t}^{\sharp\otimes\sharp}=B_{t}$.
\end{Lemma}

\begin{figure}[t]\centering
\includegraphics[clip,scale=0.4]{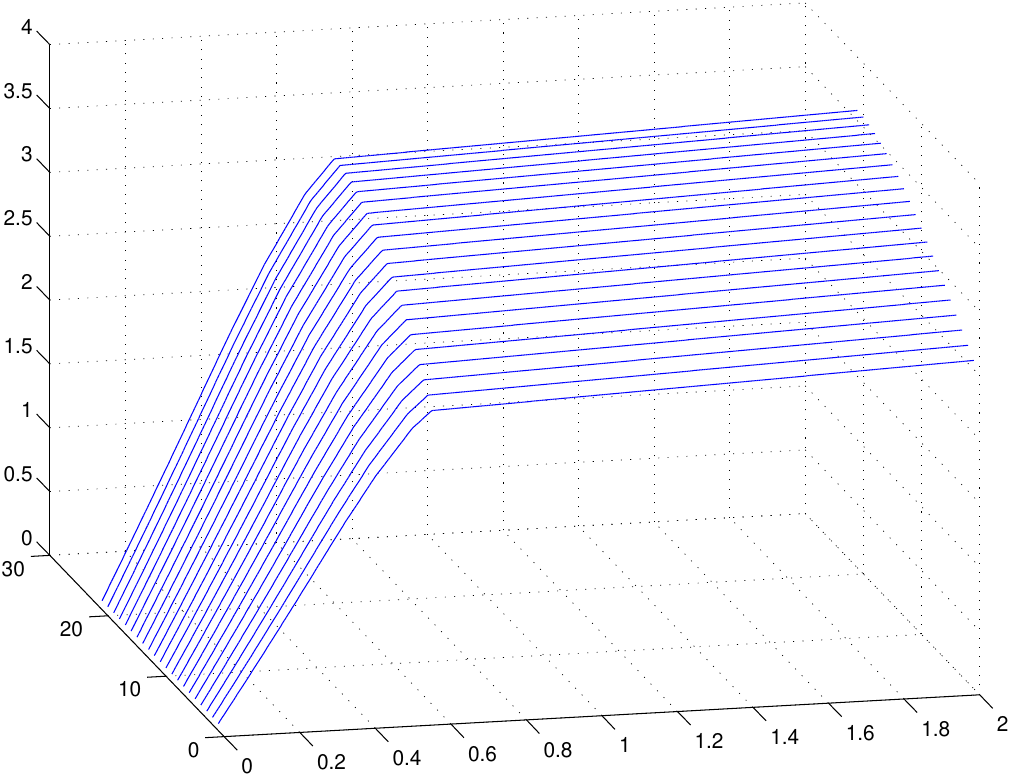}\hspace{0.3cm} \includegraphics[clip,scale=0.4]{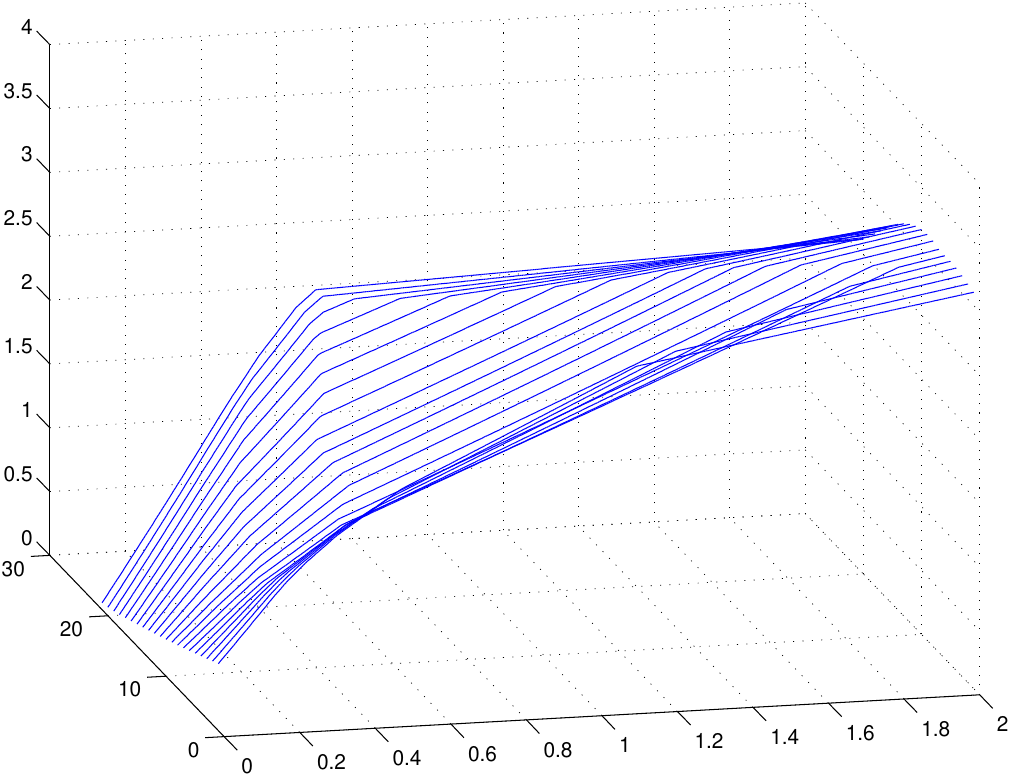}
\caption{Bounds on $\left\Vert B_{t}(U,V)^{2}-I\right\Vert $.}\label{fig:Bounds-at_t}
\end{figure}

\begin{figure}[t]\centering
\includegraphics[clip,scale=0.4]{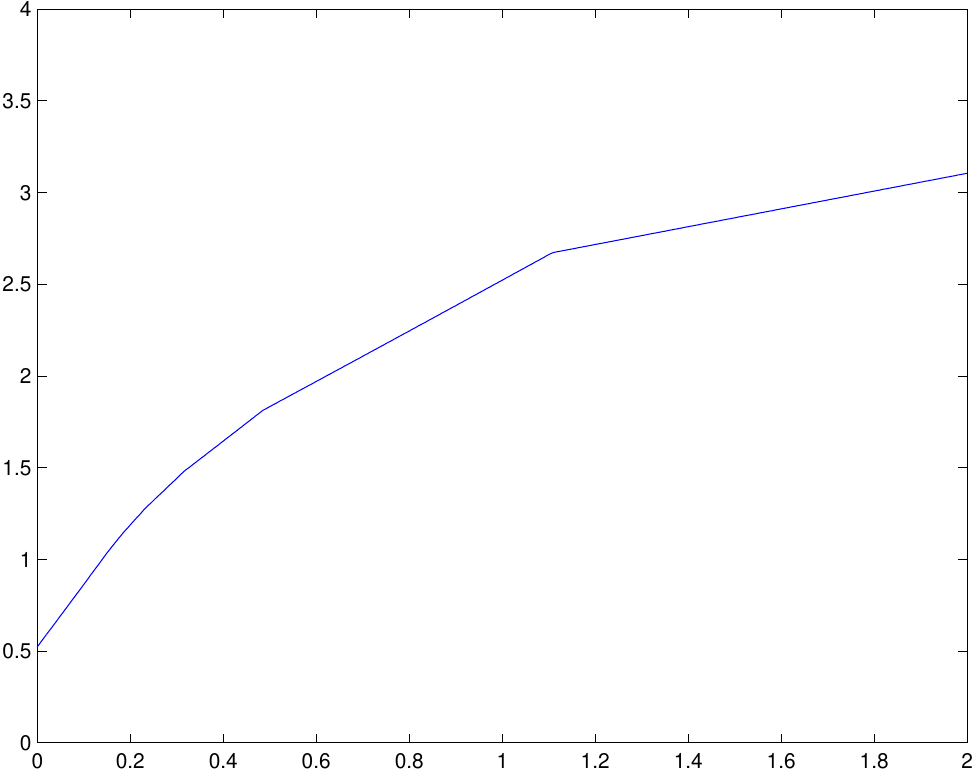}
\caption{Bounds on $\Vert B_L(U,V)^{2}-I\Vert $.}\label{fig:Bounds-at_1}
\end{figure}

\begin{proof}
We can select paths, illustrated
in Figs.~\ref{Fig6} and~\ref{Fig7}, $f_{t}$, $g_{t}$ and $h_{t}$, from the standard triple $ (f_{0},g_{0},h_{0})=(f,g,h) $, in
Def\/inition~\ref{Def:BoffByf_g_h}, to $(f_{1},g_{1},h_{1})$ where
\begin{gather}
f_{1}(x)=\frac{1}{\pi}x,\qquad g_{1}(x)=0,\qquad h_{1}(x)=\sqrt{1-\frac{1}{\pi^{2}}x^{2}}.
\label{eq:functions_for_log_method}
\end{gather}
The conditions $f^{2}+g^{2}+h^{2}=1$ and $gh=0$ hold along the path.
This gives us paths of matrices
\begin{gather*}
B_{t}(U,V)=\left[
\begin{matrix} f_{t}(K) & g_{t}(K)+\tfrac{1}{2} \{h_{t}(K),U \}
\vspace{1mm}\\
g_{t}(K)+\tfrac{1}{2} \{h_{t}(K),U^{*} \} & -f_{t}(K)
\end{matrix}
\right]
\end{gather*}
with the needed symmetries.
It remains to show these are invertible.
One needs to compute
\begin{gather}
 ( \Vert g_{t} \Vert
+1 )\eta_{h_{t}}(\delta)+\frac{1}{4} (\eta_{h_{t}}(\delta) )^{2}+\frac{1}{2}\eta_{h_{t}^{2}}(\delta)+\eta_{q_{t}}(\delta)
\label{eq:bound_at_t}
\end{gather}
and check that this takes value less than $1$ at $\delta=\frac{1}{8}$.
This is too much to do by hand, so use a~computer\footnote{See code at
\url{http://repository.unm.edu/handle/1928/23494}.} to repeatedly calculate the constants needed in
Lemma~\ref{lem:borel_f_g_h}.
We f\/ind that~\eqref{eq:bound_at_t} takes value less than $0.95$ at $\delta=\frac{1}{8}$ for all~$t$ in a~mesh
$t_{1},\dots,t_{w}$ selected so that
\begin{gather*}
 \Vert f_{t_{j}}-f_{t_{j+1}} \Vert_{\infty}+ \Vert g_{t_{j}}-g_{t_{j+1}} \Vert_{\infty}+ \Vert
h_{t_{j}}-h_{t_{j+1}} \Vert_{\infty}\leq\sqrt{1-0.95}\approx0.2236.
\end{gather*}
We can keep $f_{t}$ f\/ixed at
\begin{gather*}
f_{t}(x)=
\begin{cases}
-1, & -\pi\leq x\leq-\frac{\pi}{2},
\\
\frac{1}{128}(150\sin(x)+25\sin(3x)+3\sin(5x)), & -\frac{\pi}{2}\leq x\leq\frac{\pi}{2},
\\
1, & \frac{\pi}{2}\leq x\leq\pi,
\end{cases}
\end{gather*}
while altering $g_{t}$ from the standard~$g$ to~$0$.
The more interesting part of the path interpola\-tes~$f_{t}$ from the above to $\frac{1}{\pi}x$ while keeping $g_{t}=0$
and
\begin{gather*}
h_{t}(x)=\sqrt{1-f_{t}(x)}.
\end{gather*}
The graphs of the computed bounds are shown in Fig.~\ref{fig:Bounds-at_t}.
These bounds have been rounded up to accommodate the various errors in computing of\/fset terms when applying
Lemma~\ref{lem:eta_Bounds_m_delta_plus_b}.
The errors in computing Fourier coef\/f\/icients lead to sub-optimal results, but do not need to be accounted for as it is
the computed coef\/f\/icients that are used when applying Lemma~\ref{lem:eta_Bounds_m_delta_plus_b}.
The analysis of the error bounds is dull and omitted.
\end{proof}

It is apparent that the limitation on the constant in this result comes from the functions used in the log
method~\eqref{eq:functions_for_log_method}.
The computed bounds are shown in Fig.~\ref{fig:Bounds-at_1}.

\begin{Theorem}
%\label{thm:LogTheorem}
Suppose~$U$ and~$V$ are self-dual unitary matrices.
If $\Vert [U,V]\Vert \leq\frac{1}{8}$ then, for any~$K$ with $-\pi\leq K\leq\pi$ and $e^{iK}=V$,
\begin{gather*}
\kappa_{2}(U,V) = \mathrm{Sign}\left(\widetilde{\mathrm{Pf}}\left(\left(
\begin{matrix} \frac{1}{\pi}K & \frac{1}{2}\left\{\sqrt{I-\frac{1}{\pi^{2}}K^{2}},U\right\}
\\
\frac{1}{2}\left\{\sqrt{I-\frac{1}{\pi^{2}}K^{2}},U^{*}\right\} & -\frac{1}{\pi}K
\end{matrix}
\right) \right)\right).
\end{gather*}
\end{Theorem}

\subsection*{Acknowledgements}

The author wishes to thank Matt Hastings and Fredy Vides for discussions, both useful and entertaining.
Also he wishes to thank Robert Israel and Nick Weaver for help via MathOverf\/low.
Finally, thanks are due to the anonymous referees, whose suggestions improved the paper, especially
Sections~\ref{sec:Bott-Index} and~\ref{sec:PfaffBott}.
This work was partially supported by a~grant from the Simons Foundation (208723 to Loring).

\pdfbookmark[1]{References}{ref}

\LastPageEnding

\end{document}